\documentclass[12pt]{article}
\usepackage{amsmath}
\usepackage{amssymb}
\usepackage{cite}

\newlength{\sperr}

\newenvironment{proof}{{\settowidth{\sperr}{\it Proof}
\par\addvspace{0.3cm}\noindent\parbox[t]{1.3\sperr}
{\it Proof}\newline }}{\nopagebreak\mbox{} {\hfill $\square$}\par\addvspace{0.3cm}}

\def\a{\alpha}
\def\b{\beta}

\def\vk{\varkappa}
\def\s{\sigma}
\def\la{\lambda}
\def\om{\omega}
\def\Om{\Omega}

\def\de{\delta}
\def\t{\theta}
\def\Up{\Upsilon}
\def\vp{\varphi}
\def\ve{\varepsilon}
\def\wh{\widehat}
\def\wt{\widetilde}
\def\ov{\overline}
\def\p{\partial}

\def\BC{{\mathbb C}}
\def\BR{{\mathbb R}}

\newtheorem{Pa}{Paper}[section]
\newtheorem{Tm}[Pa]{{\bf Theorem}}
\newtheorem{La}[Pa]{{\bf Lemma}}

\newtheorem{Rk}[Pa]{{\bf Remark}}

\newtheorem{Dn}[Pa]{{\bf Definition}}
\newtheorem{Pn}[Pa]{{\bf Proposition}}

\newcommand{\CC}
{{\mathchoice {\setbox0=\hbox{$\displaystyle\rm
C$}\hbox{\hbox
to0pt{\kern0.4\wd0\vrule height0.9\ht0\hss}\box0}}
{\setbox0=\hbox{$\textstyle\rm C$}\hbox{\hbox
to0pt{\kern0.4\wd0\vrule height0.9\ht0\hss}\box0}}
{\setbox0=\hbox{$\scriptstyle\rm C$}\hbox{\hbox
to0pt{\kern0.4\wd0\vrule height0.9\ht0\hss}\box0}}
{\setbox0=\hbox{$\scriptscriptstyle\rm C$}\hbox{\hbox
to0pt{\kern0.4\wd0\vrule height0.9\ht0\hss}\box0}}}}

\title{Construction of the solution
of the inverse spectral problem for a system
depending rationally on the spectral parameter,
Borg-Marchenko-type theorem,
and sine-Gordon equation }

\author{
Alexander Sakhnovich         
\thanks{Author's work was supported by the Austrian Science Fund (FWF) under Grant  no. Y330.}
}

\date{}

\begin{document}

\maketitle

Running head: System depending rationally on spectral parameter.

\begin{abstract}
Weyl theory for a non-classical  system depending rationally  
on the spectral parameter is treated. 
Borg-Marchenko-type uniqueness theorem is proved.
The solution of the inverse problem is constructed.
An application to sine-Gordon equation
 in laboratory coordinates is given. 
\end{abstract}



\section{Introduction} \label{intro}
\setcounter{equation}{0} 

Canonical systems
\begin{equation}\label{0.1}
\frac{d}{dx}w(x,\lambda )=i\lambda JH(x)w(x,\lambda ),\quad
H \geq
0, \quad J=\left[
\begin{array}{cc}
0 & I_n \\ I_n & 0
\end{array}
\right],
\end{equation}
where $H$ are $2n \times 2n$ matrix functions and $I_n$ is the $n
\times n$ identity matrix, are classical objects of analysis,
which include  Dirac systems, matrix string equations and
Schr\"odinger equations as particular cases. For the literature on
canonical systems see, for instance, the books \cite{AD1, dB1, GoKr,
SaL3} and various references in the papers \cite{LL, Rem,
RoSa, SaA1, SaA21'}. We shall consider systems of the form
\begin{equation} \label{0.2}
y^{\prime}(x,\la) = i \sum_{k=1}^mb_k(\la-d_k)^{-1}\Big(
\beta_k(x)^* \beta_k(x)\Big)y(x,\la) , \quad b_k=\pm 1,
\quad x
\in [0,\infty),
\end{equation}
where $y^{\prime}=\frac{d}{dx}$, $\b_k(x)=\left[
\begin{array}{lr}
\b_{k1}(x)& \b_{k2}(x)
\end{array}
\right]$ are $\BC^2$-valued
differentiable vector functions such that
\begin{equation}\label{0.3}
\sup_{0<x<\infty}
\|\b_k^{\prime}(x)\|<\infty, \quad
\b_k(x)\b_k(x)^*\equiv 1, \quad 1\leq k \leq m,
\end{equation}
and $\BC$ is the complex plane.
We shall treat also a somewhat wider class of systems (\ref{0.2}),
such that the vector functions $\b_k$ satisfy relations
\begin{equation}\label{1.d12}
\sup_{0<x<l} \|\b_k^{\prime}(x)\|<\infty
\, \, {\mathrm{for}} \, {\mathrm{all}}
\, 0<l<\infty, \quad
\b_k(x)\b_k(x)^* \equiv 1, \quad 1\leq k \leq m.
\end{equation}

Systems (\ref{0.2}) generalize a subclass of canonical systems for
the important case of several poles $d_p$ with respect to the
spectral parameter $\la$. See, for instance,  interesting papers
\cite{DIKZ, Z} on systems with rational  dependence on $\la$. 
A system of the form (\ref{0.2}), where $m=2$,
can be treated as an auxiliary system for the sine-Gordon
equation in laboratory coordinates (see  Introduction in
\cite{MST} and Section \ref{Sine} here). We always assume that
\begin{equation}\label{0.4}
d_k =\ov d_k \not=d_p \quad {\mathrm{for}} \quad
k\not=p, \quad
1\leq k,p \leq m,
\end{equation}
where $\ov d_k$ is  complex conjugate to $d_k$.

The $2\times 2$ matrix function $w(x,\la)$ satisfying (\ref{0.2})
and the normalization condition, that is,
\begin{equation}\label{d1}
w^{\prime}(x,\la) = i \sum_{k=1}^mb_k(\la-d_k)^{-1}
\beta_k(x)^* \beta_k(x)w(x,\la), \quad w(0,\la) =
I_2
\end{equation}
is called the fundamental solution of (\ref{0.2}).   Different generalizations of the notion of a
Weyl function are based on the asymptotics of the fundamental
solution (see e.g. \cite{BC, BDZ, CKvA, GesST,  L, MVZ1, MVZ2, SaA11, SaA21, Y, Z}).

A $k$--th Weyl-Titchmarsh function $\phi_k(\la)$ of the system
(\ref{0.2}) was introduced in \cite{MST} on the complex disk
\[ D_M := \displaystyle{
\left\{ \la \in \BC : \left| \la - d_k - i
\frac{b_k}{M} \right| <
\frac{1}{ M} \right\}}.\] Here, it is more convinient to change
variables and use the functions $\vp_k(\mu)=\phi_k\big(\la(\mu)\big)$, where
\begin{equation}\label{1.1}
\la=d_k+\frac{b_k}{2\mu}, \quad \mu = \frac{b_k}{2(\la-d_k)}.
\end{equation}
In view of (\ref{1.1}) the inequality $\Im\mu<-M/4$ is
equivalent
to the relation $\la \in D_M$.
\begin{Dn}\label{Wf}
A function $\vp_k(\mu)$ is called a $WT_k$ function of the system
(\ref{0.2}) with the properties (\ref{1.d12})   if
and only if there exists an $M=M_k>0$ such that $\vp_k$ is holomorphic
on the half-plane $\Im\mu<-M/4$ and for all $x\in
[0,\infty)$
\begin{equation}\label{0.5}
\sup_{\Im\mu<-M/4} \left\| w\big(x,\la(\mu)\big) \left[
\begin{array}{c}\vp_k(\mu)\\ 1 \end{array}
\right]\right\|  < \infty.
\end{equation}
\end{Dn}
We shall use the notation $\zeta=\Re \mu$, $\eta= \Im
\mu$
($\mu=\zeta+i \eta$). Here $\Re$ is the real part, $\Im$ is the imaginary part, and
the real axis will be denoted by $\BR$.
When the conditions (\ref{0.3}) hold,
 the analyticity of $\vp_k$ follows automatically from (\ref{0.5}).

The direct Weyl-Titchmarsh theory for $m=2$ (two poles) was
treated and the uniqueness of the solution of the inverse problem
was proved in \cite{MST}. Here we  shall construct  this unique solution of
the inverse spectral problem ($m\geq 2$). 
Starting from the seminal work \cite{Kr} by M. Krein structured operators
were successfully used to solve inverse spectral problems. In the cases
of Krein or self-adjoint Dirac-type systems these were operators with difference
kernels \cite{AGKLS, DyI, Kr, SaA21', SaL2, SaL3}. A somewhat more
complicated structured operators will appear in this article.

An  important series of papers by F. Gesztesy, B. Simon and
coauthors on the high energy asymptotics of the Weyl functions and local
Borg-Marchenko-type
uniqueness results has initiated a growing  interest in this important domain 
(see \cite{CGR, CGZ, GesKM, GS, GZ, SaA21', Si, Si2} and references
therein). The Weyl-Titchmarsh theory for
a non-self-adjoint case (the skew-self-adjoint  Dirac type system) has been
studied
in \cite{CG, GKS, SaA11} and  the Borg-Marchenko-type results for this
system have been published in
\cite{SaA22}. Using the procedure to construct
the solution of  the inverse problem, we obtain here
a Borg-Marchenko-type theorem for another interesting
non-self-adjoint case, namely, for system (\ref{0.2}).

Finally, an application to the boundary value problem for sine-Gordon
equation in laboratory coordinates (the case of bounded solution)
will be given.

Some preliminary results on the existence and uniqueness of the 
$WT_k$ functions, representation of the fundamental solutions,
and uniqueness of  solution of the inverse problem 
are given in the next Section \ref{prel}. The structured operators,
which are necessary to construct the solution of the inverse problem,
are studied in Section \ref{S}. The construction of the
solution of the inverse problem and Borg-Marchenko-type theorem
are contained in Section \ref{Inv}. The notion of the Weyl set
is introduced in Section \ref{Set}, and the sine-Gordon equation
is treated in Section \ref{Sine}.
\section{Preliminaries} \label{prel}
\setcounter{equation}{0} To make our paper self-contained we shall
formulate in this section some results from \cite{MST}. We
formulate them for $m\geq 2$ as the proofs for the case $m\geq 2$
are  similar to the proofs for the case $m=2$ treated in
\cite{MST}. The existence and uniqueness of the Weyl functions are
stated in Theorem 2.2 \cite{MST}.
\begin{Tm} \label{WT} Let (\ref{0.2}) be a system with
coefficients $\beta_k$ which are absolutely continuous vector
functions satisfying (\ref{0.3}) and the additional condition
\begin{equation}\label{1.7}
\beta_{k1}(0)\not=0 \quad (1\leq k \leq m).
\end{equation}
Then there exist unique $WT_k$--functions $\varphi_k$ ($1\leq k
\leq m$) of the system (\ref{0.2}).
\end{Tm}
To connect $WT_k$-functions with the solutions from $L^2$, similar
to the classical Weyl functions, let $w$ be the fundamental
solution of system (\ref{0.2}). For fixed $k$ (which we sometimes
omit in the notations), we define
\begin{equation}\label{1.0}
W(x,\mu) = W_k(x,\mu) := e^{-i x \mu} Q(x) w\big(x,\la(\mu)\big) ,
\end{equation}
where $\mu$ and $\la$ are connected by the formula (\ref{1.1}) and
$Q$ is the $2\times 2$ matrix function given by
\begin{equation}\label{1.2}
Q(x) = Q_k(x) := \left[
\begin{array}{lr}
\beta_{k1}(x) \ & \beta_{k2}(x) \\ -
\ov{\beta_{k2}(x)} \ &
\ov{\beta_{k1}(x)}
\end{array}
\right] , \qquad x \in [0,\infty).
\end{equation}
Here $\beta_{kj}$ denote the entries of $\b_k$, that is,
$\b_k=[\beta_{k1} \quad \beta_{k2}]$. By (\ref{0.3}) and
(\ref{1.2}) the function $Q(x)$ is unitary-valued, that is,
\begin{equation}\label{1.3}
Q(x)^* Q(x) = Q(x) Q(x)^* = I_2, \qquad x \in [0,\infty).
\end{equation}
From  the second relation in (\ref{0.3}) and formulas (\ref{1.1})
and (\ref{1.2}) it follows that
\begin{equation}\label{1.4}
i \frac{b_k}{\la-d_k} \, Q_k(x) \beta_k(x)^* \beta_k(x)
Q_k(x)^* =
2 i \mu \left[
\begin{array}{cc}1&0\\0&0\end{array}
\right], \quad \la=\la(\mu).
\end{equation}
By (\ref{d1}), (\ref{1.0}), and (\ref{1.4}) the matrix function
$W$ satisfies the system
\begin{equation}\label{1.5}
W^{\prime}(x,\mu) = \left( i \mu j + \xi(x,\mu)
\right) W(x,\mu),
\quad j=\left[
\begin{array}{cc}
1 & 0 \\ 0 & -1
\end{array}
\right], \quad x \in [0,\infty),
\end{equation}
where
\begin{equation}\label{1.6}
\xi(x,\mu)=Q_k^{\prime}(x)Q_k(x)^*+ iQ_k(x)\left(\sum_{p
\not=
k}\frac{b_p\b_p(x)^*\b_p(x)}{\la -d_p} \right)Q_k(x)^*.
\end{equation}
From  (\ref{1.0}),  (\ref{1.3}), (\ref{1.5}),  and (\ref{1.6})
it
follows that
\begin{equation}\label{1.d0}
W(x,\ov \mu)^*W(x, \mu)=Q(0)^*Q(0)=I_2, \quad W(0,
\mu)=Q(0).
\end{equation}
The analog of Theorem 2.4 \cite{MST}, which is formulated below,
states that $W_k(x,\mu)\left[
\begin{array}{c}
\vp_k(\mu) \\ 1
\end{array}
\right] \in L^2_2$.
\begin{Tm} \label{L} Let (\ref{0.2}) be a system with
coefficients $\beta_k$ which are absolutely continuous vector
functions satisfying (\ref{0.3}) and (\ref{1.7}). Then the
$WT_k$--functions $\varphi_k$ are unique functions such that for
some $M_k>0$ and  all $\mu$ satisfying inequality
$\Im\mu<-M_k$ we
have
\begin{equation}\label{1.9}
\int_0^{\infty}[\ov \vp_k(\mu) \quad
1]W_k(x,\mu)^*W_k(x,\mu)\left[
\begin{array}{c}
\vp_k(\mu) \\ 1
\end{array}
\right]dx<\infty.
\end{equation}
\end{Tm}
We shall need some details from the proof of Theorem 2.2
\cite{MST} (Theorem \ref{WT} here) in our further considerations.
Notice that the Dirac-type system can be written down in the form
(\ref{1.5}), where $\xi$ does not depend on $\mu$. Similar to the
Dirac-type system case \cite{SaA11}, choose a value $M>0$ such
that
\begin{equation}\label{1.8}
\sup_{x\in [0,\infty), \Im\mu<-M/4}
\|\xi(x,\mu)\|<\frac{1}{4}M.
\end{equation}
  By (\ref{1.5}) and (\ref{1.8})
one can see that for $\Im\mu<-M/4$ we have \cite{MST}:
\begin{equation}\label{1.10}
\frac{d}{d x}R(x,\mu)>0, \quad
R(x,\mu):=Q(0)W(x,
\mu)^*jW(x,\mu)Q(0)^*.
\end{equation}
Now, put
\begin{equation}\label{1.d1}
{\mathfrak A}(x,\mu)={\mathfrak A}_k(x,\mu)=\{ {\mathfrak A}_{jp}(x,\mu)
\}_{j,p=1}^2:=
Q(0)W_k(x, \mu)^{-1}.
\end{equation}
According to the second relation in (\ref{1.d0}) and to
(\ref{1.10}) we have $R(x,\mu)>j$ or, equivalently, ${\mathfrak
A}(x,\mu)^*j{\mathfrak A}(x,\mu)<j$. Thus, the linear fractional
transformation
\begin{equation}\label{1.d2}
\psi_k(l,\mu)=\frac{{\mathfrak
A}_{11}(l,\mu)\t(\mu)+ {\mathfrak
A}_{12}(l,\mu)}{{\mathfrak A}_{21}(l,\mu)
\t(\mu)+{\mathfrak
A}_{22}(l,\mu)}, \quad |\t(\mu)|\leq 1, \quad l>0,
\quad \Im \mu
<-\frac{M}{4},
\end{equation}
where $\t$ is a holomorphic parameter function, is well-defined, and
\begin{equation}\label{1.d3}
|\psi_k(l,\mu)|<1.
\end{equation}
The class of functions $\psi_k(l,\cdot)$ given by (\ref{1.d2}) is
denoted by ${\cal N}_k(l)$. Using (\ref{1.10}), it is shown in
\cite{MST} that ${\cal N}_k(l_1) \subset {\cal N}_k(l_2)$ for
$l_1>l_2$.  For each $l>0$ and for each $\mu$ ($ \Im \mu <
-
\frac{ M}{4}$) the values of $ \psi_k(l,\mu)$ ($\psi_k \in 
{\cal
N}_k(l)$) can be parametrized
\begin{equation}\label{1.d4}
\psi_k(l,\mu) = \rho_1(l,\mu)^{-1/2} \Theta(l,\mu)
\rho_2(l,\mu)^{-1/2} + \rho_0(l,\mu), \qquad l \in
(0,\infty),
\end{equation}
\begin{equation}\label{1.d5}
\rho_0 = -R_{11}^{-1} R_{12}, \quad \rho_1 = R_{11}, \quad
\rho_2
= \left( R_{21} R_{11}^{-1} R_{12} - R_{22} \right)^{-1},
\end{equation}
where $|\Theta(l,\mu)|\leq 1$ and $R_{jp}$ are the entries of
$R=\{R_{jp}\}_{j,p=1}^2$.
The set of values of $ \psi_k(l,\mu)$ ($\psi_k \in
{\cal N}_k(l)$) coincides with the disk on the right-hand side of
(\ref{1.d4}), that is, the values of $\psi_k$ form the  so called
Weyl disks. The functions $\rho_1(l)^{-1/2}$ and $\rho_2(l)^{-1/2}
$ are decreasing, and for $\rho_1$ we have
\begin{equation}\label{1.d7}
\rho_1(l)\geq 1-2l\Big(\frac{M}{4}+ \Im \mu\Big)
\to \infty, \quad {\mathrm{when}} \quad
l \to \infty.
\end{equation}
Therefore, the  intersection of the Weyl disks in (\ref{1.d4}) is
a Weyl point, that is, there is only one function $\wt \psi_k(\mu)$,
which belongs to all ${\cal N}_k(l)$:
\begin{equation}\label{1.d6}
\bigcap_{l<\infty}{\cal N}_k(l)=: \wt \psi_k(\cdot).
\end{equation}
To study the asymtotics of $\wt \psi_k$ we need the representation
of the fundamental solution $W$ from Theorem 2.1 \cite{MST}:
\begin{Tm} \label{Repr}
Let $\beta_k(x)$ be absolutely continuous  $\BC^2$-valued vector functions
on the interval $[0, \, l]$ $( 0<l<\infty)$ satisfying relations 
\begin{equation}\label{1.d18}
\sup_{0<x<l} \|\b_k^{\prime}(x)\|<\infty ,
\quad
\b_k(x)\b_k(x)^* \equiv 1, \quad 1\leq k \leq m.
\end{equation}
Then $W(x,\mu)$ $(x \in [0, \, l])$ of the form (\ref{1.0})
admits
a representation
\begin{equation}\label{1.d13}
\begin{array}{l}
W_k(x,\mu)Q(0)^* = \displaystyle{ e^{i \mu x j} \left( D_k(x) +
\sum_{p\not=k}\left( \mu - \frac{b_k}{2(d_p-d_k)}
\right)^{-1}
D_p(x) \right)}  \\ \displaystyle{ + \int_{-x}^x e^{i
\mu u}
N_k(x,u) d u +\sum_{p\not=k} \left( \mu -
\frac{b_k}{2(d_p-d_k)}
\right)^{-1} \int_{-x}^x e^{i \mu u} N_p(x,u) d u} \\
\displaystyle{+ {\rm
O}(\mu^{-2}),}
\end{array}
\end{equation}
for $\mu=\zeta+i\eta,  \eta\not= 0, \ |\zeta|
\to \infty$, where
$D_s$ $(1\leq s \leq m)$ are continuous diagonal matrix functions,
$D_k^* = D_k^{-1}$, and
\begin{equation}\label{1.d14}
\sup_{ |u|\leq x\leq l} \left(  \sum_{s=1}^m\|
N_s(x,u)\| \right)
< \infty.
\end{equation}
\end{Tm}
From the representation (\ref{1.d13}), after some calculations one
gets
\begin{equation}\label{1.d15}
\lim_{\Im \mu \to - \infty}\wt \psi_k(\mu)=0
\end{equation}
uniformly with respect to $\Re\mu$ (see formula (2.29) in
\cite{MST}). As (\ref{1.d15}) holds and
$\b_{k1}(0)\not=0$, there is a value $\wt M$ such that
\begin{equation}\label{1.d8}
|\ov{\b_{k2}(0)}\wt
\psi_k(\mu)+\b_{k1}(0)|>\ve>0, \quad \Im
\mu<-\wt M \leq -\frac{M}{4}.
\end{equation}
The $WT_k$-function is defined in the domain $\Im \mu<-\wt M$ by
the formula
\begin{equation}\label{1.d10}
\vp_k(\mu)=\frac{\ov{\b_{k1}(0)}\wt
\psi_k(\mu)-\b_{k2}(0)}{\ov{\b_{k2}(0)}\wt
\psi_k(\mu)+\b_{k1}(0)}.
\end{equation}
Finally, using (\ref{1.d4}) the estimate
\begin{equation}\label{1.d11}
|\psi_k(l,\mu)-\breve \psi_k(l,\mu)|\leq
2\exp\Big(\big(i(\ov \mu
- \mu)+M/2 \big)l\Big), \quad \Im
\mu<-\frac{M}{4}
\end{equation}
was proved for  the arbitrary $\psi_k, \breve \psi_k \in
{\cal
N}_k(l)$ (see Lemma 2.3 \cite{MST}).  This estimate, in its turn,
helped to prove that the solution of the inverse problem is unique
(Theorem 3.1 \cite{MST}):
\begin{Tm} \label{Uniq}
For given $WT_k$-functions $\vp_k$ ($1\leq k\leq m$) there is at
most one system (\ref{0.2}) satisfying conditions (\ref{0.3}) and
(\ref{1.7}).
\end{Tm}

In this paper we shall construct also the solution of the following inverse
problem.
\begin{Dn} \label{DnInv}
The inverse spectral problem  for system (\ref{0.2}), (\ref{1.d12}) 
is the problem to recover  the system,
that is, to recover the matrix function 
\[\b(x)=\left[
\begin{array}{c}
\b_1(x) \\ \cdots \\ \b_m(x)
\end{array}
\right],
\]
such that the relations (\ref{1.d12})  and 
(\ref{1.7}) hold and that the given functions $\vp_k$ are system's
$WT_k$-functions. We call $\b$ the potential of  system (\ref{0.2}).
\end{Dn}
The corresponding uniqueness theorem  
 was proved in \cite{MST} using Theorem \ref{Repr}.
\begin{Tm} \label{Uniq2}
For given functions $\vp_k$ $(1\leq k\leq m)$,
which admit asymptotic representations
\begin{equation}\label{1.d16}
\vp_k(\mu)=c_k+O(\mu^{-1}), 
\quad c_k\in \BC,
\quad \Im
\mu<-\frac{M}{4}, \quad \mu \to \infty,
\end{equation}
there is at
most one system (\ref{0.2}) satisfying conditions (\ref{1.d12})  and 
(\ref{1.7})  and such that the functions $\vp_k$ are the system's
$WT_k$-functions. 
Moreover, it follows from (\ref{1.d16}) that
\begin{equation}\label{1.d17}
c_k=-\b_{k2}(0)/\b_{k1}(0).
\end{equation}
\end{Tm}

\section{S-nodes} \label{S}
\setcounter{equation}{0} Interesting developments of the classical
results on the inverse problems were obtained using the notion of
the $S$-node (see \cite{SaL2, SaL3} and references therein). The
non-self-adjoint systems were also recovered from their Weyl
functions using $S$-nodes \cite{SaA11, SaA21, SaA22}. Here, we
shall consider  $S$-nodes corresponding to system (\ref{0.2}). By
$\{{\cal H}_1,\, {\cal H}_2\}$ we denote the class of the
linear
bounded operators acting from the Hilbert space ${\cal H}_1$ into
the Hilbert space ${\cal H}_2$, diag means diagonal matrix, and
$A_0(l)\in \{L^2_m(0,l), \, L^2_m(0,l)\}$ is an integration
operator: $(A_0f)(x)=\int_0^x f(u) du$, $x\in [0,l]$. As in the
expression $A_0f$ above, sometimes we omit $l$ in our notations.
Later we shall also denote by $A_0(l)$ the operator of integration
in $L^2(0,l)$. We denote the identity operator by $I$. It is always 
clear from the context, where $I$ is acting.

Now, we introduce three bounded operators $A(l)$, $S(l)$, and $\Pi(l)$.
Operators $A(l)$ have simple structure and do not depend on the
choice of $\b_k(x)$:
\begin{equation}\label{2.1}
Af=A(l)f=Df+iBA_0(l)f, \quad A(l) \in \{L^2_m(0,l), \,
L^2_m(0,l)\},
\end{equation}
\begin{equation}\label{2.2}
D={\rm{diag}}\{d_1,d_2,\ldots,d_m\}, \quad
B={\rm{diag}}\{b_1,b_2,\ldots,b_m\}.
\end{equation}
Operator $\Pi$ is an operator of multiplication by the $m \times
2$ matrix function $[\Phi_1(x) \quad \Phi_2(x)]$:
\begin{equation}\label{2.3}
\Pi(l) \left[
\begin{array}{c}
g_1 \\ g_2
\end{array}
\right]=g_1\Phi_1(x)+g_2\Phi_2(x), \quad \Pi(l) \in
\{\BC^2, \,
L^2_m(0,l)\},
\end{equation}
where the entries of the absolutely continuous column vector
functions $\Phi_p$ ($p=1,2$) are denoted by $\Phi_{kp}$, and we
require
\begin{equation}\label{2.4}
\Phi_{k1}(x)\equiv 1, \quad  \Phi_{k2}^{\prime}(x)\in
L^2(0,l),
\quad  1 \leq k \leq m.
\end{equation}
Later we shall recover $\Phi_2$ from the Weyl functions. 

Operator
$S \in \{L^2_m(0,l), \, L^2_m(0,l)\}$ is chosen so that it
satisfies the operator identity
\begin{equation}\label{2.5}
AS-SA^*=i\Pi\Pi^*.
\end{equation}
Therefore we say that $A$, $S$, and $\Pi$ form an $S$-node.
\begin{Pn}\label{Pn2.1}
Let operators $A$ and $\Pi$ be given by equalities
(\ref{2.1})-(\ref{2.3}), and let $\Pi$ satisfy (\ref{2.4}).
Then
the unique bounded operator $S$, which satisfies (\ref{2.5}), has
the form
\begin{equation}\label{2.6}
Sf=B\wt D f+\int_0^ls(x,u)f(u)du, \quad
s(x,u)=\{s_{kp}(x,u)\}_{k,p=1}^m,
\end{equation}
where
\begin{equation}\label{2.6'}
\wt D
=I_m+{\rm{diag}}\{|\Phi_{12}(0)|^2,|\Phi_{22}(0)|^2,\ldots,|\Phi_{m2}(0)|^2\}.
\end{equation}
The entries of $s(x,u)$ on the main diagonal are defined by the
equalities
\begin{equation}\label{2.7}\displaystyle
s_{kk}(x,u)=\frac{b_k}{2}\int_{|x-u|}^{x+u}\Phi_{k2}^{\prime}\left(
\frac{v+x-u}{2}\right)\overline{\Phi_{k2}^{\prime}\left(
\frac{v+u-x}{2}\right)}dv.
\end{equation}
The offdiagonal $(k\not=p)$ entries of $s(x,u)$ are defined by the
equalities
\begin{equation}\label{2.d2}
s_{kp}(x,u)=\frac{1}{2\pi}\int_{\Gamma}\Big((\la
I-A_k)^{-1}[1
\quad \Phi_{k2}]\Big)(x)\Big( (\ov \la I-A_p)^{-1}[1
\quad
\Phi_{p2}]\Big)(u)^*d\la,
\end{equation}
where $\Gamma=\{\la: |\la -d_k|=\varepsilon>0\}$ is
anticlockwise
oriented, $\varepsilon<|d_k-d_p|$, $I$ is the identity operator,
$A_k=d_kI+ib_kA_0$, and we use $A_0$ here to denote the
integration operator in $L^2(0,l)$.
\end{Pn}

\begin{proof}
Write down $S$ in the matrix form 
\[
S=\{S_{kp}\}_{k,p=1}^m, \qquad
S_{kp}\in \{L^2(0,l), \, L^2(0,l)\}. 
\]
Then identity
(\ref{2.5})
takes the form
\begin{equation}\label{2.5'}
A_0 S_{kk}+S_{kk}A_0^*=b_k\Pi_k\Pi_k^*, \quad  \Pi_k \left[
\begin{array}{c}
g_1 \\ g_2
\end{array}
\right]=g_1+g_2\Phi_{k2}(x), \quad \Pi_k \in
\{\BC^2, \,
L^2(0,l)\},
\end{equation}
\begin{equation}\label{2.5''}
A_kS_{kp}-S_{kp}A_p^*=i\Pi_k\Pi_p^*, \quad
A_k=d_kI+ib_kA_0,\quad
k \not=p.
\end{equation}
The bounded solution $T$ of the equation $TA_0+A_0^*T=Q$ for $Q$
of the form $Q=\int_0^lq(x,t)\cdot dt$ is constructed in Theorem
1.3 (p.11) \cite{SaL2'}.  After easy transformations we derive
from this result the solution of (\ref{2.5'}) too. Namely, we have
\begin{align}\label{2.d3}
&S_{kk}f= \frac{b_k}{2} 
 \frac{d}{d x}\int_0^l\frac{\p}{\p u}\left(
\int^{x+u}_{|x-u|}\Big(
1+\Up(v,x,u) \Big)d
v\right)f(u)du, \\
&
\Up(v,x,u):=\Phi_{k2}\big((v+x-u)/2\big)\ov{\Phi_{k2}\big((v-x+u)/2\big)}.
\nonumber
\end{align}
From (\ref{2.d3}), taking into account the second relation in
(\ref{2.4}), we derive (\ref{2.6})-(\ref{2.7}). To get
(\ref{2.d2}) rewrite (\ref{2.5''}) in the form
\begin{equation}\label{2.0}
(\la I-A_k)^{-1}S_{kp}-S_{kp}(\la I-A_p^*)^{-1}=i(\la
I-A_k)^{-1}\Pi_k\Pi_p^*(\la I-A_p^*)^{-1},
\end{equation}
and notice that $\s(A_k)=d_k$ and $\s(A_p^*)=d_p$, where
$d_k\not=d_p$ and $\s$ means spectrum. So, we recover  $S_{kp}$ by
integration of the both parts of (\ref{2.0}) in the small
neighborhood of $d_k$:
\begin{equation}\label{2.0'}
S_{kp}=\frac{1}{2\pi}\int_{\Gamma}(\la
I-A_k)^{-1}\Pi_k\Pi_p^*(\la
I-A_p^*)^{-1}d\la .
\end{equation}
In other words we have
\begin{equation}\label{2.d1}
S_{kp}=\int_0^ls_{kp}(x,u)\cdot du,
\end{equation}
where $s_{kp}$ satisfies (\ref{2.d2}).
\end{proof}
\begin{Rk}\label{res}
It is easy to check the explicit formula for the resolvent of
$A_k$:
\begin{eqnarray}\label{2.9}
\Big((\la I -A_k)^{-1}f\Big)(x) \hspace{22em}
\\
=(\la -d_k)^{-1}f(x)+ib_k(\la -d_k)^{-2}\int_0^x\exp
\left(\frac{ib_k(x-u)}{\la -d_k}\right)f(u)du. \nonumber
\end{eqnarray}
\end{Rk}
Denote by $P_r$ ($r \leq l$) the orthogonal projector from
$L^2_m(0,l)$ onto $L^2_m(0,r)$, that is, let $P_r \in \{L^2_m(0,l),
\, L^2_m(0,r)\}$ and let $(P_rf)(x)=f(x)$ for $x\in (0,r)$. Notice
that $P_rA(l)=A(r)P_r$. Therefore, we get
\begin{equation}\label{2.11}
A(r)P_rSP_r^*-P_rSP_r^*A(r)^*=i\Pi(r)\Pi(r)^*,
\end{equation}
by applying $P_r$ from the left and $P_r^*$ from the right  to the
both parts of (\ref{2.5}).
\begin{Rk}\label{red}
According to (\ref{2.11}), the unique operator  $S(r)$ satisfying
the identity
\begin{equation}\label{2.12}
A(r)S(r)-S(r)A(r)^*=i\Pi(r)\Pi(r)^*, \quad r<l,
\end{equation}
is given by the formula
\begin{equation}\label{2.13}
S(r)=P_rSP_r^*=B\wt D +\int_0^r s(x,u) \cdot du,
\end{equation}
where $s$ does not depend on $r$ and $B\wt D$ means the operator
of multiplication by the  matrix $B\wt D$.
\end{Rk}
We shall need some properties of $S(l)$.
\begin{Pn}\label{Pn2.2}
The operator $S$ constructed in Proposition \ref{Pn2.1} is
self--adjoint, boundedly invertible and $S^{-1}$ admits a triangular
factorization
\begin{equation}\label{2.10}
S^{-1}=V^*BV, \quad (Vf)(x)=\wt
D^{-\frac{1}{2}}f(x)+\int_0^xV(x,u)f(u)du,
\end{equation}
where
\begin{equation}\label{2.14}
V(r,u)=B\wt D^{\frac{1}{2}}T_r(r,u) \quad (r\geq u),
\end{equation}
and $T_r$ is the matrix kernel of the integral operator
\begin{equation}\label{2.15}
T(r)=S(r)^{-1}=B\wt D^{-1}+\int_0^rT_r(x,u) \cdot du.
\end{equation}
\end{Pn}

\begin{proof}
The operator $S$ is self-adjoint as the unique solution of
(\ref{2.5}). (One could also prove it by
(\ref{2.6})-(\ref{2.d2}).)
The invertibility of $S$ is proved by contradiction. Suppose that
$S$ is not invertible. In view of  the special structure
(\ref{2.6}) of $S$, it means that $S$ has an eigenvector
$f\not=0$, such that $Sf=0$.  Taking into account identity
(\ref{2.5}) and equality $Sf=0$, we derive $\big(f,
\Pi\Pi^*f\big)_{L^2}=0$, where $\big(\cdot, \cdot
\big)_{L^2}$
denotes the scalar product in $L^2_m(0,l)$. It is immediate that
$\Pi^* f=0$. Apply the both parts of (\ref{2.5}) to $f$ and use
the equalities $Sf=0$, $\Pi^*f=0$ to obtain $SA^*f=0$. So, from
Sf=0 it follows that $SA^*f=0$. In other words, we have SL=0 for
the linear span $L$ of the vectors $(A^*)^kf$ ($k \geq 0$).
Therefore, we have $\dim \, L<\infty$. As $A^*L \subseteq L$
and
$\dim \, L<\infty$, there is an eigenvector $g$ of $A^*$:
$A^*g=cg$, $g\not=0$, and $g \in L$. Hence, by the definition of
$A_k$ in (\ref{2.5''}), there is an eigenvector of integration in
$L^2(0,l)$: $A_0^*g_k=\wt c g_k$, $g_k\not=0$. This is impossible,
and so we come to a contradiction, that is, $S$ is invertible.

In view of (\ref{2.12}) and (\ref{2.13}) the invertibility of the
operators $P_rSP_r^*$ ($r<l$) is proved quite similar to the
invertibility of $S$. By (\ref{2.7}) and (\ref{2.d2}) the function
$s(x,u)$ is continuous. Thus, the factorization conditions from
"result 2" Section IV.7 \cite{GoKr} are fulfilled for $\wt
S^{-1}$, where $\wt S=B\wt D^{-\frac{1}{2}}S\wt
D^{-\frac{1}{2}}$.
Hence, the factorization formula for $S^{-1}$ in (\ref{2.10}), the
second relation in  (\ref{2.10}), and equality (\ref{2.14}) follow.
\end{proof}
\begin{Rk} \label{Rk3.5} Let the conditions of Proposition
\ref{Pn2.1}
be fulfilled and put
\begin{equation}\label{2.d6}
\b(x)=\left[
\begin{array}{c}
\b_1(x) \\ \cdots \\ \b_m(x)
\end{array}
\right]=(V\Phi)(x), \quad \Phi(x):=[\Phi_1(x) \quad
\Phi_2(x)] \quad (0 \leq x \leq l),
\end{equation}
where $V$ is applied to $\Phi(x)$  columnwise. In other words, we have 
\begin{equation}\label{2.d6'}
V\Pi g=\b(x)g \quad (g\in \BC^2).
\end{equation}
Then the matrix functions $\b_k$ satisfy  the second relation in
(\ref{1.d12}),
that is, $\b_k\b_k^* \equiv 1$. Indeed,
from (\ref{2.5}) and the first equality in (\ref{2.10}), it follows
that
\[
V^*BVA-A^*V^*BV=iV^*BV\Pi\Pi^*V^*BV,  \quad {\mathrm{i.e.}}, 
\]
\begin{equation}\label{2.d10}
VAV^{-1}B-B(V^*)^{-1}A^*V^*=iV\Pi\Pi^*V^*.
\end{equation}
By the definitions (\ref{2.1}) and (\ref{2.10}) of $A$ and $V$, the
operator
$VAV^{-1}B$ is lower triangular and has the form $DB+\int_0^x
\gamma(x,u)\cdot du$.
The operator $B(V^*)^{-1}A^*V^*$  
is upper triangular. Hence, one can derive the kernel $\gamma$ of the
integral term of
$VAV^{-1}B$
from (\ref{2.d6'}) and (\ref{2.d10}).
We get
\begin{equation}\label{2.d11}
(VAV^{-1}f)(x)=Df(x)+i\b(x)\int_0^x\b(u)^*Bf(u)du.
\end{equation}
By (\ref{2.d11}) it is immediate that
\begin{equation}\label{2.d11'}
VA=DV+i\b(x)\int_0^x\b(u)^*BV \cdot du.
\end{equation}
Rewrite (\ref{2.d11'}) as the equality of the kernels of the
corresponding integral operators:
\[
i\wt D^{-\frac{1}{2}}B +V(x,u)D+i\int_u^x V(x,v) B dv
\hspace{15em}
\]
\begin{equation}\label{2.d12}
=DV(x,u)+i\b(x)\Big(\b(u)^*B\wt D^{-\frac{1}{2}}+
\int_u^x \b(v)^*BV(v,u)  dv\Big). 
\end{equation}
When $x=u$, the equality of the main diagonals of the both sides of
(\ref{2.d12}) implies
$\b_k\b_k^* \equiv 1$.
\end{Rk}
Now, introduce the transfer matrix functions in the Lev Sakhnovich
form \cite{SaL1, SaL2, SaL3}:
\begin{equation}\label{2.d4}
w_A(r,\la)=I_2-i\Pi(r)^*S(r)^{-1}\big(A(r)-\la I)^{-1}\Pi(r).
\end{equation}
The following lemma is essential for the solution of the inverse problem.
\begin{La}\label{Law}
Let the conditions of Proposition \ref{Pn2.1} be fulfilled, and
let the $S$-node be given by the formulas (\ref{2.1})-(\ref{2.4})
and (\ref{2.6})-(\ref{2.d2}). Then we have
\begin{equation}\label{2.d5}
\frac{d}{dr}w_A(r,\la)=i\b(r)^*B(\la I_m
-D)^{-1}\b(r)w_A(r,\la).
\end{equation}
\end{La}
\begin{proof}  First, introduce several notations. Let $P_1(r,\delta)$
and
$P_2(r,\delta)$ denote ortoprojectors from $L^2_m(0,r+\delta)$ on
$L^2_m(0,r)$ and $L^2_m(r,r+\delta)$, respectively. That is, let
$P_1(r, \delta)f \in L^2_m(0,r)$, $P_2(r, \delta)f \in
L^2_m(r,r+\delta)$, and
\begin{equation}\label{2.d9}
  \Big(P_1(r, \delta)f\Big)(x)=f(x), \,\, 0<x<r;
\quad \Big(P_2(r, \delta)f\Big)(x)=f(x), \, \,
r<x<r+\delta .
\end{equation}
For operators $K$ acting in $L^2_m(0,r+\delta )$  we put
$K_{jp}:=P_j(r, \delta)KP_p(r, \delta)^*$ ($j,p=1,2$). In
particular, we use notations
\begin{equation}\label{2.d8}
T(u):=S(u)^{-1}, \quad T_{22}:=P_2(r, \delta)T(r+\delta)P_2(r,
\delta)^*.
\end{equation}
According to \cite{SaL1} (see also Theorem 2.1 from Chapter 1 in \cite{SaL3})
we have
\begin{eqnarray}\label{2.d7}
w_A(r+\delta,\la)&-&w_A(r,\la) =-i\Pi(r+\delta
)^*S(r+\delta )^{-1}
\\ &&\times (A_{22}-\la I)^{-1}T_{22}^{-1}P_2(r,
\delta)S(r+\delta
)^{-1}\Pi(r+\delta ) w_A(r,\la). \nonumber
\end{eqnarray}
Using formula (\ref{2.d6}) and the first equality in (\ref{2.10}),  we
rewrite (\ref{2.d7}) in the form
\begin{equation}\label{2.d13}
w_A(r+\delta,\la)-w_A(r,\la)=i\int_r^{r+\delta}\b(x)^*(Z_\delta\b)(x)
dx \, w_A(r,\la),
\end{equation}
where the operator $Z_\delta \in \{L^2_m(r,r+\delta),
\, L^2_m(r,r+\delta)\}$ is given by the formula
\begin{equation}\label{2.d14}
Z_\delta=BV(\la I-A_{22})^{-1}T_{22}^{-1}P_2(r, \delta)V^*B.
\end{equation}
It is easy to see that $Z_\delta-B(\la I-D)^{-1}$ is an integral
operator, and
we shall show below that the kernel of this operator is bounded:
\begin{equation}\label{2.d15}
Z_\delta=B(\la I-D)^{-1}+\int_r^{r+\delta}z_\delta(x,u)
\cdot du, \quad \sup\|z_\delta (x,u)\|<\infty,
\end{equation}
where $r\leq x,u \leq r+\delta\leq l$.
For that purpose  notice that the kernel $s(x,u)$ is continuous,
and according to Section IV.7 \cite{GoKr} the kernel $T_r(x,u)$ is
continuous with respect to $r$, $x$, and $u$ ($x,u\leq r\leq l$)
too. Hence, the functions $V(x,u)$ and $\b(x)$ are continuous, and
we also have
\begin{equation}\label{2.d16}
\sup_{x,u \leq l}\|s(x,u)\|<\infty, \quad
\sup_{x,u\leq r\leq l}\|T_r(x,u)\|<\infty.
\end{equation}
From the definition (\ref{2.d8}) we get
\begin{equation}\label{2.d17}
T_{22}^{-1}=S_{22}-S_{21}S_{11}^{-1}S_{12},
\end{equation}
and so, by (\ref{2.d16}) the kernel of the integral term of $T_{22}^{-1}$
is bounded. In view of (\ref{2.9}), for the $k$-th entry of $(\la
I-A_{22})^{-1}f$ we have
\begin{eqnarray}\label{2.d18}
\Big((\la I-A_{22})^{-1}f\Big)_k(x) \hspace{22em}
\\
=(\la -d_k)^{-1}f(x)+ib_k(\la -d_k)^{-2}\int_r^x\exp
\left(\frac{ib_k(x-u)}{\la -d_k}\right)f(u)du, \nonumber
\end{eqnarray}
and the kernel of $(\la I-A_{22})^{-1}-(\la I-D)^{-1}$ is bounded for any
fixed $\la$ ($\la \not=d_k$) too.
Therefore, for $Z_\delta$ given by (\ref{2.d14}) the formula 
(\ref{2.d15}) is true. Recall that $\b$ is continuous.
Thus, from (\ref{2.d13}) and (\ref{2.d15}) it follows that
\begin{equation}\label{2.d19}
\lim_{\delta \to +0}\delta^{-1}\Big(
w_A(r+\delta,\la)-w_A(r,\la) \Big)=i\b(r)^*B(\la I_m
-D)^{-1}\b(r)w_A(r,\la).
\end{equation}
Quite similar one can prove that
\begin{equation}\label{2.d20}
\lim_{\delta \to +0}\delta^{-1}\Big( w_A(r,\la)-
w_A(r-\delta,\la) \Big)=i\b(r)^*B(\la I_m
-D)^{-1}\b(r)w_A(r,\la).
\end{equation}
By (\ref{2.d19}) and (\ref{2.d20}) equality (\ref{2.d5}) holds.
\end{proof}
\begin{Rk}\label{Rkw} It is easy to see that definition
(\ref{2.d4})
implies $\lim_{r \to 0}w_A(r,\la)=I_2$. Hence, the matrix function
$w_A$, which is treated in Lemma \ref{Law}, is the fundamental solution
of the system (\ref{0.2}) corresponding to $\b(x)=(V\Phi)(x)$.
\end{Rk}
\section{Inverse problems: \\ 
construction of the solution} \label{Inv}
\setcounter{equation}{0}
According to formulas (\ref{1.d3}), (\ref{1.d8}), and
(\ref{1.d10}), one can choose a sufficiently large value $M>0$,
so that all the $WT_k$-functions ($1\leq k \leq m$) are well-defined 
and bounded in  the half-plane $\Im\mu<-M/4$. 
\begin{Dn} \label{Dn41}
The bounded
matrix function
\begin{equation}\label{2.16}
\vp(\mu)={\mathrm{col}}[\vp_1(\mu), \quad
\vp_2(\mu), \quad
\ldots, \quad \vp_m(\mu)],
\end{equation}
\begin{equation}\label{2.19}
\sup_{\Im \mu < -M/4 }\| \vp(\mu)
\|<\infty,
\end{equation}
where col means column,  is called the Weyl function of system
(\ref{0.2}).  
\end{Dn}
One of our main results  is the next theorem.
\begin{Tm}\label{TmInv1}
System (\ref{0.2}) satisfying conditions (\ref{0.3}) and
(\ref{1.7}) is uniquely recovered from its Weyl function. To
recover (\ref{0.2}) on an arbitrary fixed interval $[0, \, l]$ we
use  the following procedure.

First, introduce the column vector function $\Phi_2(x)$
$(0<x<\infty)$ by the Fourier transform:
\begin{equation}\label{2.17}
\Phi_2(x)=\frac{i}{2
\pi}\int_{-\infty}^{\infty}\mu^{-1}e^{i \mu
x}\vp\left(\frac{\mu}{2}\right)d\zeta \quad
(\mu=\zeta +i \eta,
\quad \eta<-\frac{M}{2}).
\end{equation}
That is, we define $\Phi_2$ on $(0, \, l)$ and $(0,\, \infty)$
via
norm limit l.i.m. in $L_2^m$:
\begin{equation}\label{2.17'}
\displaystyle{\Phi_2(x)=\frac{i}{2 \pi}e^{-\eta
x}{\mathrm{l.i.m.}}_{a \to \infty}
\int_{-a}^{a}\mu^{-1}e^{i \zeta
x} \vp\left(\frac{\mu}{2}\right)d\zeta.}
\end{equation}
Here the right-hand side of (\ref{2.17'}) equals $0$ for $x<0$.

Next, substitute $\Phi_2(x)$ into formulas
(\ref{2.6'})-(\ref{2.d2}) to introduce operator $S\in
\{L^2_m(0,l), \, L^2_m(0,l)\}$ of the form (\ref{2.6}). 
Apply formulas (\ref{2.10})-(\ref{2.15}) to recover operator $V$
from the operators $S(r)=P_rSP_r^*$ ($r \leq l$). Then, the matrix
functions
$\b_k(x)$ and therefore, system (\ref{0.2}) is recovered  by the
formula
(\ref{2.d6}), that is,
$\b(x)=V[\Phi_1(x) \quad \Phi_2(x)]$,
where $\Phi_1$ is given by  the
first relation in (\ref{2.4}). 
\end{Tm}
\begin{proof} 
{\bf Step 1.} Let system (\ref{0.2}) satisfy conditions (\ref{0.3})
and
(\ref{1.7}) and let $\vp$ be the  Weyl function of system
(\ref{0.2}).
According to the inequality (\ref{2.19}),  the vector function
$\Phi_2$ is well-defined by  (\ref{2.17}) and does not depend on
the choice of $\eta<-\frac{M}{4}$. So, we can fix some
$\eta<-\frac{M}{4}$. To show that $\Phi_2$ is absolutely
continuous, introduce functions
\begin{equation}\label{2.22}
\wh \psi_k (\mu):=\frac{{\mathfrak
A}_{12}(l,\mu)}{{\mathfrak
A}_{22}(l,\mu)}\in{\cal N}_k(l), \quad \wh
\vp_k(\mu)=\frac{\ov{\b_{k1}(0)}\wh
\psi_k(\mu)-\b_{k2}(0)}{\ov{\b_{k2}(0)}\wh
\psi_k(\mu)+\b_{k1}(0)}.
\end{equation}
By (\ref{1.d0}) and (\ref{1.d1}) we have ${\mathfrak
A}(l,\mu)=Q_k(0)W_k(l, \ov \mu)^*$. Hence, by Theorem
\ref{Repr} the
matrix function $\wh \psi_k $ admits representation
\begin{equation}\label{2.20}
\wh \psi_k (\mu)=\big(\ov
{d_{k22}(l)}\big)^{-1}\int_{-l}^le^{-i\mu
(l+u)}\ov{N_{k21}(l,u)}du+g_1(\mu), \end{equation} where $d_{k22}$
is the corresponding entry of the diagonal unitary matrix $D_k$,
$N_{k21}$ is an entry of the bounded matrix function $N_k$, and
the function $g_1(\zeta +i\eta) \in L^1(-\infty, \,
\infty)\cap
L^{\infty}(-\infty, \, \infty)$ with respect to the variable
$\zeta$. (Recall that $L^{\infty}$ is the space of bounded
functions.) Taking into account (\ref{1.d8}) and (\ref{1.d11}),
without loss of generality we can assume
\begin{equation}\label{2.23}
|\ov{\b_{k2}(0)}\wh
\psi_k(\mu)+\b_{k1}(0)|>\ve>0, \quad \Im
\mu<-\frac{M}{4}.
\end{equation}
In view of relations (\ref{2.22})-(\ref{2.23}) we get
\begin{equation}\label{2.24}
\wh
\vp_k(\mu)=-\frac{\b_{k2}(0)}{\b_{k1}(0)}+\frac{|\b_{k1}(0)|^2-|\b_{k2}(0)|^2}{\b_{k1}(0)^2}\wh
\psi_k (\mu)+g_2(\mu),
\end{equation}
where $g_2(\zeta +i\eta) \in L^1(-\infty, \,
\infty)\cap
L^{\infty}(-\infty, \, \infty)$. Finally, notice that uniformly
for $x$ on the intervals $[\delta, \, l]$ we have
\begin{equation}\label{2.25}
\lim_{a \to \infty}\frac{i}{2
\pi}\int_{-a}^{a}\mu^{-1}e^{i \mu
x}d\zeta =-1 \quad (\mu=\zeta +i \eta, \quad
\eta<-\frac{M}{2}).
\end{equation}
From (\ref{2.20}), (\ref{2.24}), and (\ref{2.25}) it follows that
the function
\begin{equation}\label{2.26}
\wh \Phi_{k2}(x):=\frac{i}{2
\pi}\int_{-\infty}^{\infty}\mu^{-1}e^{i
\mu x}\wh \vp_k\left(\frac{\mu}{2}\right)d\zeta
\quad (\mu=\zeta
+i \eta, \quad \eta<-\frac{M}{2}),
\end{equation}
is absolutely  continuous and
\begin{equation}\label{2.36}
\sup_{0<x<l}|\wh \Phi_{k2}^{\prime}(x)|<\infty.
\end{equation}
Next, we shall show that
\begin{equation}\label{2.26'}
\wh \Phi_{k2} (x)= \Phi_{k2}(x)  \quad {\mathrm{for}}
\quad 0\leq x
\leq l,
\end{equation}
where $\Phi_{k2}$ denotes the $k$-th entry of the $\BC^m$-valued vector function
$\Phi_2$. The proof of (\ref{2.26'}) requires some considerations.
According to the definitions (\ref{1.d10}) and (\ref{2.22}) of
$\vp_k$ and $\wh \vp_k$, respectively, we have
\begin{equation}\label{2.27}
\vp_k(\mu) - \wh
\vp_k(\mu)=\frac{\big(|\b_{k1}(0)|^2+|\b_{k2}(0)|^2\big)\big(\wt
\psi_k(\mu) - \wh \psi_k(\mu)\big)}
{\big(\ov{\b_{k2}(0)}\wt
\psi_k(\mu)+\b_{k1}(0)\big)\big(\ov{\b_{k2}(0)}\wh
\psi_k(\mu)+\b_{k1}(0)\big)}.
\end{equation}
In view of (\ref{1.d8}), (\ref{1.d11}), and (\ref{2.23}), formula
(\ref{2.27}) implies for some $C, \wt M>0$ that
\begin{equation}\label{2.28}
|\vp_k(\mu) - \wh \vp_k(\mu)|<C(l)e^{2\eta l},
\quad \eta<-\wt
M(l)<-\frac{M}{4} \quad (\mu=\zeta + i \eta).
\end{equation}
Taking into account (\ref{2.19}) and (\ref{2.28}), one can see
that functions $\mu^{-1}
\vp_k\left(\frac{\mu}{2}\right)$ and
$\mu^{-1}\wh \vp_k\left(\frac{\mu}{2}\right)$ are
holomorphic  in
the half-plane $\eta= \Im \mu<-2\wt M(l)$ and belong to
$L^2(-
\infty, \, \infty)$ with respect to $\zeta=\Re \mu$ for
any fixed
$\eta<-2\wt M(l)$. Therefore, according to Theorem V \cite{PW}
these functions admit Fourier representations
\begin{equation}\label{2.29}
\mu^{-1}
\vp_k\left(\frac{\mu}{2}\right)={\mathrm{l.i.m.}}_{a
\to
\infty}\int_0^ae^{-i \mu x}e^{2x\wt M}f(x)dx, \quad f
\in L^2(0,
\infty),
\end{equation}
\begin{equation}\label{2.30}
\mu^{-1} \wh
\vp_k\left(\frac{\mu}{2}\right)={\mathrm{l.i.m.}}_{a
\to \infty}\int_0^ae^{-i \mu x}e^{2x\wt M}\wh f(x)dx,
\quad \wh f
\in L^2(0, \infty).
\end{equation}
Using Plansherel's theorem and formulas (\ref{2.17})
and(\ref{2.26}), we express $f(x)$ and $\wh f(x)$ in (\ref{2.29})
and (\ref{2.30}) via $\Phi_{k2}(x)$ and $\wh \Phi_{k2}(x)$,
respectively, and obtain:
\begin{equation}\label{2.31}
\frac{i}{\mu}
\vp_k\left(\frac{\mu}{2}\right)=\int_0^{\infty}e^{-i
\mu x}\Phi_{k2}(x)dx, \quad e^{-2x\wt M}\Phi_{k2}(x)
\in L^2(0,
\infty),
\end{equation}
\begin{equation}\label{2.32}
\frac{i}{\mu} \wh
\vp_k\left(\frac{\mu}{2}\right)=\int_0^{\infty}e^{-i
\mu x}\wh
\Phi_{k2}(x)dx, \quad e^{-2x\wt M}\wh \Phi_{k2}(x) \in
L^2(0,
\infty).
\end{equation}
Consider now the entire matrix function
\begin{equation}\label{2.33}
Y_0(\mu)=e^{i \mu l}\int_0^{l}e^{-i \mu
x}\Big(\Phi_{k2}(x)-\wh
\Phi_{k2}(x)\Big)dx =\int_0^{l}e^{i \mu
(l-x)}\Big(\Phi_{k2}(x)-\wh
\Phi_{k2}(x)\Big)dx.
\end{equation}
From (\ref{2.31})-(\ref{2.33}) it follows that for $\mu=\zeta +
i
\eta$, $\eta <-2\wt M$ we have
\begin{equation}\label{2.34}
|Y_0(\mu)|\leq \int_l^{\infty}e^{\eta
(x-l))}\Big(|\Phi_{k2}(x)|+|\wh \Phi_{k2}(x)|\Big)dx+
e^{-\eta
l}\left|\frac{i}{\mu}\Big(
\vp_k\left(\frac{\mu}{2}\right)- \wh
\vp_k\left(\frac{\mu}{2}\right)\Big)\right|.
\end{equation}
Recall that $\Phi_{k2}, \, \wh \Phi_{k2} \, \in \,
L^2(0, \infty)$.
Hence, from (\ref{2.28}) and (\ref{2.34}) we derive
\begin{equation}\label{2.35}
\sup_{ \eta<-2\wt M-\ve}|Y_0(\mu)|< \infty
\quad (\ve>0), \quad
\lim_{\eta \to -\infty}|Y_0(\mu)|=0 \quad
(\mu=\zeta + i \eta).
\end{equation}
By the definition (\ref{2.33}) and by the first relation in
(\ref{2.35}), the entire functon $Y_0$ is bounded in $\BC$. So, in
view of the second relation in (\ref{2.35}) we have $Y_0=0$, and
equality (\ref{2.26'}) is immediate. 

From (\ref{2.36}) and
(\ref{2.26'}) we get the second relation in (\ref{2.4}).
Therefore, $\Pi$ satisfies the conditions of Proposition
\ref{Pn2.1}, and the $S$-node and the operator $V$ in our theorem  are
well-defined.
As the conditions of Propositions \ref{Pn2.1} and \ref{Pn2.2}
are satisfied, the matrix function $\breve \b(x)=(V\Phi)(x)$ is
well-defined and unique.
To prove the theorem means to prove the equalities
\begin{equation}\label{0p}
\breve \b_k(x)^* \breve \b_k(x) \equiv 
 \b_k(x)^*  \b_k(x) \quad (1\leq k \leq m).
\end{equation}

{\bf Step 2.}
By  (\ref{2.26}), taking into account (\ref{2.20}), (\ref{2.24}),
and (\ref{2.25}),
we have
\begin{equation}\label{2.37}
e^{\eta x}\wh \Phi_{k2}^{\prime}(x)=-\frac{1}{2
\pi}{\mathrm{l.i.m.}}_{a \to \infty}\int_{-a}^{a}e^{i
\zeta x}\left(\wh
\vp_k\Big(\frac{\mu}{2}\Big)+\frac{\b_{k2}(0)}{\b_{k1}(0)}\right)d\zeta,
\end{equation}
\begin{equation}\label{2.38}
e^{\eta x}\wh \Phi_{k2}^{\prime}(x) \in L^2(0, \infty)
\cap L^{\infty}(0, \infty), \quad \eta<-\frac{M}{2}.
\end{equation}
Introduce a $2 \times 2$ matrix function:
\begin{equation}\label{1p}
 \Om(\mu)=[\Om_1(\mu) \quad
\Om_2(\mu)]=Q(0)^*\left[
\begin{array}{cc}
\exp\{-2ir\mu\} &\wh \psi_k(\mu) \\ 0
& 1
\end{array}
\right].
\end{equation}
First, let us show that for  $r \leq l$ and for sufficiently large $M$ the
inequalities
\begin{equation}\label{2p}
\sup_{\Im \mu
<-M/4}\|w_A\big(r,\la(\mu)\big)\Om_1(\mu)\|<\infty , \quad  
\sup_{\Im \mu
<-M/4}\|w_A\big(r,\la(\mu)\big)\Om_2(\mu)\|<\infty
\end{equation}
are true. Indeed,  by  Lemma \ref{Law} $w_A$ satisfies 
system (\ref{0.2}), where we substitute $\breve \b_k$ instead
of $\b_k$. Hence,
for  sufficiently large $M$ we get
\begin{eqnarray}\label{3p}
&&\frac{d}{dr}\Big((\exp\{-2ir(\mu-\ov
\mu)\})w_A(r,\la)^*w_A(r,\la)\Big)
\\ && \nonumber
=\exp\{-2ir(\mu-\ov \mu)\}w_A(r,\la)^*
\Big(4\eta\big(
I_2-\breve\b_k(r)^*\breve\b_k(r)\big)+g_3(r,\la)\Big)w_A(r,\la),
\\ && \nonumber
\la=\la(\mu),\quad g_3(r,\la)<CI_2 \quad (\eta=\Im \mu <-M/4),
\end{eqnarray}
for some $C>0$, which does not depend on $r$. As $\eta
(I_2-\breve\b_k(r)^*\breve\b_k(r))\leq 0$,
formula  (\ref{3p}) implies
\begin{equation}\label{4p}
\sup_{\Im \mu <-M/4}\|\exp\{-2ir\mu\}
w_A(r,\la)\|<\infty ,
\end{equation}
that is, the first inequality in  (\ref{2p}) holds.
Next, let us prove  that
\begin{equation}\label{y1}
\sup_{\Im \mu <-M/4}\|w_A(r,\la)
\left[
\begin{array}{c}
\wh \vp_k(\mu) \\  1
\end{array}
\right]
\|<\infty \quad (r\leq l).
\end{equation}
For that purpose consider $w_A$ again and notice that from  the operator
identity (\ref{2.12})
and definition (\ref{2.d4}) it follows \cite{SaL2} that
\begin{eqnarray}\label{n.20}
&&w_A(r,\la)^*w_A(r,\la) \\
&&=I_2-i(\la -\ov \la)\Pi(r)^*\big(A(r)^*-\ov
\la I\big)^{-1}S(r)^{-1}
(A(r)- \la I\big)^{-1}\Pi(r). \nonumber
\end{eqnarray}
By   (\ref{2.4}) and  (\ref{2.9})
 we derive the equality: 
\begin{equation}\label{n.21}
\big(A_k(r)-\la I\big)^{-1}[\Phi_{k1} \quad
\Phi_{k2}]=-2b_k\mu e^{2i\mu x}\Big[1 \quad 
\Phi_{k2}(0)+\int_0^x e^{-2i \mu
u}\Phi_{k2}^{\prime}(u)du\Big],
\end{equation}
which helps to estimate   the right-hand side
of   (\ref{n.20}).
According to (\ref{2.32}) and  (\ref{2.38}) the equality
\begin{equation}\label{y2}
\wh
\vp_k\left({\mu}\right)=-\wh\Phi_{k2}(0)-\int_0^{\infty}e^{-2i
\mu u}\wh \Phi_{k2}^{\prime}(u)du
\end{equation}
is true.  In view of (\ref{2.26'}), (\ref{n.21}), and
(\ref{y2}) we have
\begin{equation}\label{y3}
\big(A_k(r)-\la I\big)^{-1}\Pi_k\left[
\begin{array}{c}\wh \vp_k(\mu)\\ 1 \end{array}
\right]=2b_k\mu e^{2i\mu x}\int_x^\infty e^{-2i \mu
u}\wh \Phi_{k2}^{\prime}(u)du.
\end{equation}
By (\ref{2.38}) and  (\ref{y3}) for sufficiently large $M$ we get
\begin{equation}\label{y4}
\left\|\big(A_k(r)-\la I\big)^{-1}\Pi_k\left[
\begin{array}{c}\wh \vp_k(\mu)\\ 1 \end{array}
\right]\right\|\leq C(r)|\mu/\eta|.
\end{equation}
Moreover, for sufficiently large $M$ the functions $\wh \vp_k(\mu)$
and resolvents
$\big(A_p(r)-\la I\big)^{-1}$ ($p\not=k$) are uniformly bounded
in the domain $\Im \mu<-M/4$.
Therefore, the inequalities (\ref{y4})
imply
\begin{equation}\label{y5}
\left\|\big(A(r)-\la I\big)^{-1}\Pi\left[
\begin{array}{c}\wh \vp_k(\mu)\\ 1 \end{array}
\right]\right\|\leq C_1(r)|\mu/\eta|.
\end{equation}
Notice that for $\la=d_k+\frac{b_k}{2\mu}$ we have $|\la -\ov
\la |=-\eta|\mu|^{-2}$. 
Hence, from  (\ref{n.20}) and (\ref{y5}) follows (\ref{y1}).  

By (\ref{2.22}) we have
\begin{equation}\label{y5'}
\left[
\begin{array}{c}\wh \vp_k(\mu)\\ 1 \end{array}
\right]=Q(0)^*\left[
\begin{array}{c}\wh \psi_k(\mu)\\ 1 \end{array}
\right]\big(\ov{\b_{k2}(0)}\wh
\psi_k(\mu)+\b_{k1}(0)\big)^{-1}.
\end{equation}
Substitute (\ref{y5'}) into (\ref{y1}) and use (\ref{1.d3}) to
derive the second relation
in (\ref{2p}). Thus, (\ref{2p}) is valid.

Now, let us show that
\begin{equation}\label{y6}
\sup_{\Im \mu <-M/4}\|w\big(r,\la(\mu)\big)\Om(\mu)
\|<\infty.
\end{equation}
Here, the inequality
\begin{equation}\label{y7}
\sup_{\Im \mu <-M/4}\|\exp\{-2ir\mu\}
w(r,\la)\|<\infty , \quad \la=\la(\mu)
\end{equation}
is proved similar to  (\ref{4p}).  

As $\vp$ of the form  (\ref{2.16}) is the Weyl function,
so $\vp_k$ satisfies  (\ref{0.5}).  The inequality
\begin{equation}\label{y8}
\sup_{\Im \mu <-M/4}\|w(r,\la)
\left[
\begin{array}{c}
\wh \vp_k(\mu) \\  1
\end{array}
\right]
\|<\infty \quad (r \leq l)
\end{equation}
follows for sufficiently large $M$ from (\ref{0.5}), (\ref{2.28}), and
(\ref{y7}).
In view of (\ref{1.d3}) and (\ref{y5'}),  formula (\ref{y8})
yields 
\begin{equation}\label{y9}
\sup_{\Im \mu
<-M/4}\|w(r,\la)\Om_2(\mu)\|<\infty.
\end{equation}
By  (\ref{y7}) and (\ref{y9})  inequality  (\ref{y6}) is valid.

{\bf Step 3.} Here,  using (\ref{2p}) and (\ref{y6})  we shall show
that
\begin{equation}\label{y9'}
\sup_{|\Im \mu|
>M/4}\|w_A(r,\la)w(r,\la)^{-1}\|<\infty.
\end{equation}
First, taking into account  (\ref{1.0}) and (\ref{1.d13}), we obtain
\begin{equation}\label{y10}
w(r,
\la)Q(0)^*=e^{ir\mu}Q(r)^*W(r,\mu)Q(0)^*=Q(r)^*D_k(r){\mathrm{diag}}\{
e^{2ir\mu},1\}+g_4(\mu), 
\end{equation}
where the $2 \times 2$ matrix function $g_4(\zeta +i\eta)$ belongs
$L^2_{2\times 2}(-\infty, \, \infty)$ for each fixed
$\eta<-M/4$, that is, the entries of
$g_4$ belong $L^2$.
 According to
 (\ref{1.d3}) and (\ref{2.20}) the function $\wh \psi(\zeta
+i\eta)$, where
 $\eta<-M/4$, is bounded
 and belongs to $L^2(-\infty, \, \infty)$ with respect to
$\zeta$.
 Therefore, formulas (\ref{1p}) and (\ref{y10}) imply that for
$\eta <-M/4$ we have
 \begin{equation}\label{y11}
w(r,\la)\Om(\mu)=Q(r)^*D_k(r)+g_5(\mu), \quad g_5(\zeta
+i\eta)\in L^2_{2\times 2}(-\infty, \, \infty)
.
\end{equation}
 After some evident change of
variables in  Theorem VIII from \cite{PW}, 
by  (\ref{y6}) and (\ref{y11}) we can apply it
to $w(r,\la)\Om(\mu)-Q(r)^*D_k(r)$. That is, we  derive for
$\eta<-M/4$
the representation
 \begin{equation}\label{y12}
w(r,\la)\Om(\mu)=Q(r)^*D_k(r)+\int_0^{\infty}e^{-ix\mu}f(x)dx,
\quad e^{\eta x}f(x)\in L^2_{2\times 2}(-\infty, \,
\infty). 
\end{equation}
In view of  (\ref{2p}) and (\ref{y12}) for sufficiently large $M$ we
have
\begin{equation}\label{y13}
\sup_{\Im \mu
<-M/4}\|w_A(r,\la)w(r,\la)^{-1}\|<\infty.
\end{equation}
According to (\ref{d1}) and Lemma \ref{Law}  the equalities
\begin{equation}\label{y14}
w(r,\ov \la)^*=w(r,\la)^{-1}, \quad w_A(r,\ov
\la)^*=w_A(r,\la)^{-1}
\end{equation}
hold. In particular, we have
\begin{equation}\label{y14'}
\frac{1}{\det w(r, \ov \la)^*}=\det w(r,\la), \quad
\frac{1}{\det w_A(r, \ov \la)^*}=\det w_A(r,\la).
\end{equation}
Recall also that $w$ and $w_A$ are $2 \times 2$ matrices, and so 
(\ref{y14}) and
 (\ref{y14'}) yield
\begin{equation}\label{y15}
w(r, \la)=\big(w(r, \ov \la)^*\big)^{-1}
=\big(\det w(r,\la)\big)jJ\ov{w(r, \ov \la)}Jj
\end{equation}
\begin{equation}\label{y16}
w_A(r, \la)=\big(\det w_A(r,\la)\big)jJ\ov{w_A(r,
\ov \la)}Jj.
\end{equation}
Here we put $n=1$ in the definition of $J$, that is,
\begin{equation}\label{y16'}
 J=\left[
\begin{array}{cc}
0 & 1\\ 1 & 0
\end{array}
\right], \quad  j=\left[
\begin{array}{cc}
1 & 0 \\ 0 & -1
\end{array}
\right].
\end{equation}
From  (\ref{0.3}) and (\ref{d1})  it follows that
\begin{eqnarray} \nonumber
&&\det w(r,\la)=\exp\left(
i\int_0^r{\mathrm{Tr}}\Big(\sum_{k=1}^mb_k(\la-d_k)^{-1}
\beta_k(x)^* \beta_k(x)\Big)dx\right)
\\  \label{y17}
&&=\exp\Big( ir\sum_{k=1}^mb_k(\la-d_k)^{-1}
\Big),
\end{eqnarray}
where Tr means trend. In a similar way from Lemma \ref{Law} and Remark
\ref{Rk3.5}
we get
\begin{equation}\label{y18}
\det w_A(r,\la)=\exp\Big(
ir\sum_{k=1}^mb_k(\la-d_k)^{-1}\Big).
\end{equation}
Finally, substitute (\ref{y15})-(\ref{y18}) into (\ref{y13}) to
obtain
 \begin{equation}\label{y19}
\sup_{\Im \mu <-M/4}\|w_A(r,\ov \la)w(r,\ov
\la)^{-1}\|<\infty.
\end{equation}
Using (\ref{1.1}) one can see that inequalities  (\ref{y13}) and
(\ref{y19}) imply
 (\ref{y9'}).

{\bf Step 4.} 
Taking into account (\ref{d1}),  (\ref{1.1}),  (\ref{y14}), and
Lemma \ref{Law},
it is easy to see that for
\begin{equation}\label{y20}
K(r, \mu):= w_A\big(r, \la(\mu)\big)w\big(r, \la(\mu)\big)^{-1}
\end{equation}
we have
\begin{equation}\label{y21}
\underline{\lim}_{R \to \infty} \frac{1}{\ln R} \ln
\ln \sup_{|\mu|=R}\|K(r, \mu)\|\leq 1.
\end{equation}
In view of  (\ref{y9'}) and (\ref{y21}), we can apply to
$K(r,\mu)$ the Phragmen-Lindel\"of  theorem.
Thus, we see that $\|K(r,\mu)\|$ is uniformly bounded for
sufficiently large $|\zeta|$ on the strip
$|\eta|\leq M/4$. Using this and inequality  (\ref{y9'}),  we have
\begin{equation}\label{y22}
\sup_{|\mu|>M_1}\|K(r, \mu)\|< \infty
\end{equation}
for some $M_1>0$. 

Let us switch to the variable $\la$. From (\ref{1.1}),  (\ref{y20}), and 
(\ref{y22}) we derive
that 
\[
\wt K(r,\la):= w_A(r, \la)w(r, \la)^{-1}
\]
 is bounded in the
deleted neighborhood of
$\la =d_k$. That is, by  Theorem 11.4 from \cite{Silv} the matrix function
$\wt K(r,\la)$ is holomorphic
at $\la =d_k$. Since $k$ ($1\leq k\leq m$) is arbitrary,
$\wt K(r,\la)$ is an entire function with 
respect to the variable $\la$. Moreover, it easy to see that
\begin{equation}\label{y23}
\lim_{\la \to \infty} w_A(r, \la)=\lim_{\la \to
\infty} w(r, \la)=\lim_{\la \to \infty} w_A(r, \la)w(r,
\la)^{-1} = I_2.
\end{equation}
By Liouville's theorem, it follows from  (\ref{y23}) that the entire
function
$w_Aw^{-1}$ is constant, that is, 
\begin{equation}\label{y24}
w_A(r, \la)  \equiv w(r, \la) \quad (0\leq r \leq l).
\end{equation}
Identity   (\ref{y24}) implies (\ref{0p}).
\end{proof}
The construction of the solution of the inverse problem, which is described in
Definition \ref{DnInv},
is given below.
\begin{Tm}\label{TmInv2}
Let a  vector function $\vp(\mu)$ be holomorphic in the half-plane
$\Im \mu<-\frac{M}{4}$ and admit there an asymptotic
representation
\begin{equation}\label{n.1}
\vp(\mu)=\a_0+\frac{\a_1}{\mu}+O\Big(\frac{1}{\mu^{2}}\Big),
\quad \Im \mu<-\frac{M}{4}, \quad \mu \to
\infty.
\end{equation}
Then $\vp$ is a Weyl function of the unique
system (\ref{0.2}) satisfying conditions (\ref{1.d12}) and 
(\ref{1.7}).

To
recover (\ref{0.2}) on an arbitrary fixed interval $[0, \, l]$ one can
use  the  same procedure as in Theorem \ref{TmInv1}.
First, introduce the column vector function $\Phi_2(x)$
($0<x<\infty$) by the Fourier transform (\ref{2.17}) or,
equivalently, (\ref{2.17'}).
Here the right-hand side of (\ref{2.17'}) equals $0$ for $x<0$.

Next, substitute $\Phi_2(x)$ into formulas
(\ref{2.6'})-(\ref{2.d2}) to introduce operator $S\in
\{L^2_m(0,l), \, L^2_m(0,l)\}$ of the form (\ref{2.6}). 
Apply formulas (\ref{2.10})-(\ref{2.15}) to recover operator $V$
from the operators $S(r)=P_rSP_r^*$ ($r \leq l$). 

Then, the matrix
functions
$\b_k(x)$ and therefore, system (\ref{0.2}) is recovered  using the
formula
(\ref{2.d6}), that is,
$\b(x)=V[\Phi_1(x) \quad \Phi_2(x)]$,
where $\Phi_1$ is given by  the
first relation in (\ref{2.4}). 
\end{Tm}
\begin{proof} 
{\bf Step 1.} Let us show that the potential $\b=V[\Phi_1 \quad
\Phi_2]$ 
is differentiable and relations (\ref{1.d12}) and (\ref{1.7}) hold.

According to (\ref{2.25})  and to the asymptotic relation (\ref{n.1}), 
the vector function
$\Phi_2$ is well-defined by  (\ref{2.17}) and does not depend on
the choice of $\eta<-\frac{M}{2}$. So, we can fix some
$\eta<-\frac{M}{2}$. Moreover, one can easily see that
$\Phi_2(x)$ is twice differentiable and $e^{\eta
x}\Phi_2^{\prime \prime}(x) \in L^2_m(0, \infty)$.
Therefore, the matrix function $s(x,u)$ given by (\ref{2.7}) and
(\ref{2.d2})
is continuous. So, $T_r(x,u)$ given by (\ref{2.15}) is continuous with
respect to $r$, $x$, and $u$ \cite{GoKr}. 

The product of the right-hand
sides of the first equality in (\ref{2.6}) (for
$l=r$) and of the formula (\ref{2.15}) equals $I_m$, which
can be written down as the equality
\begin{equation}\label{n.2}
s(x,u)B\wt D^{-1}+B\wt DT_r(x,u)+\int_0^rs(x,t)T_r(t,u)dt=0
\quad (x,u \leq r)
\end{equation}
for the kernels of the integral operators. In view of the continuity
of $s$ and $T_r$, equality (\ref{n.2}) is true pointwise. Changing the order
of multiplication of the right-hand sides of (\ref{2.6}) (for
$l=r$) and (\ref{2.15}) we get also
\begin{equation}\label{n.3}
B\wt D^{-1}s(x,u)+T_r(x,u)B\wt D+\int_0^rT_r(x,t)s(t,u)dt=0
\quad (x,u \leq r).
\end{equation}
From (\ref{2.10}), (\ref{2.14}), and (\ref{2.d6}) it follows that
\begin{equation}\label{n.4}
\b(r)=\wt D^{-\frac{1}{2}}[\Phi_1(r) \quad
\Phi_2(r)]+B\wt D^{\frac{1}{2}}\int_0^rT_r(r,u)
[\Phi_1(u) \quad \Phi_2(u)]du.
\end{equation}
Using (\ref{n.2})-(\ref{n.4}), one can show that
\begin{equation}\label{n.5}
\b^{\prime}(r)=\wt D^{-\frac{1}{2}}[0 \quad
\Phi_2^{\prime}(r)]
+B\wt D^{\frac{1}{2}}\big(Y_1+Y_2+Y_3\big),
\end{equation}
where
\begin{eqnarray}\label{n.6}
&&Y_1=T_r(r,r)\Phi(r), \quad \Phi(r)=[\Phi_1(r)
\quad \Phi_2(r)], \\
&&Y_2=-T_r(r,r)
\int_0^r\Big(S(r)^{-1}s(x,r)\Big)^*\Phi(x)dx,\label{n7}
\\
&& 
Y_3=-B\wt
D^{-1}\int_0^r\Big(S(r)^{-1}\frac{\partial}{\partial
u}s(x,u)\big|_{u=r}\Big)^*
\Phi(x)dx, \label{n.8}
\end{eqnarray}
and $S(r)^{-1}$ is applied to the matrix functions columnwise.

Indeed, let us prove that 
\begin{equation}\label{n.9}
\frac{d}{dr}\int_0^rT_r(r,u)
\Phi(u)du=Y_1+Y_2+Y_3.
\end{equation}
It is immediate that for $\de>0$ we have
\begin{eqnarray}\nonumber
&&\int_0^{r+\de}T_{r+\de}(r+\de,u)
\Phi(u)du-\int_0^{r}T_{r}(r,u)
\Phi(u)du
\\ \nonumber &&
=\int_r^{r+\de}T_{r+\de}(r+\de,u)
\Phi(u)du
+\int_0^{r}\Big(T_{r+\de}(r+\de,u)-T_{r+\de}(r,u)\Big)\Phi(u)du
\\&&
+
\int_0^{r}\Big(T_{r+\de}(r,u)-T_{r}(r,u)\Big)\Phi(u)du.\label{n.10}
\end{eqnarray}
As $T_r(x,u)$ is continuous, we get
\begin{equation}\label{n.11}
\lim_{\de \to
0}\de^{-1}\int_r^{r+\de}T_{r+\de}(r+\de,u)
\Phi(u)du=Y_1.
\end{equation}
Substitute $r+\de$ instead of $r$ into (\ref{n.2}) to obtain the
equality
\[
s(x,u)B\wt D^{-1}+B\wt
DT_{r+\de}(x,u)+\int_0^{r+\de}s(x,t)T_{r+\de}(t,u)dt=0,
\]
which can be rewritten as
\begin{equation}\label{n.12}
S(r)T_{r+\de}(x,u)=-s(x,u)B\wt
D^{-1}-\int_r^{r+\de}s(x,t)T_{r+\de}(t,u)dt
\end{equation}
for all fixed values of $u$ ($u \leq r$).
According to (\ref{n.12}) we have
\begin{eqnarray}\nonumber
&&T_{r+\de}(x,
r+\de)-T_{r+\de}(x,r)=-S(r)^{-1}\Big(\big(s(x,r+\de)-s(x,r)\big)B\wt
D^{-1}\\
&&\label{n.13}
+\int_r^{r+\de}s(x,t)
\Big(T_{r+\de}(t,r+\de)-T_{r+\de}(t,r)\Big)dt\Big).
\end{eqnarray}
Recall that $\Phi$ is twice differentiable. Hence, for $u \geq r$
and $x \leq r$ the matrix function $s(x,u)$ is differentiable with
respect to $u$ and $ \frac{\partial}{\partial u}s(x,u)$ is
continuous with respect to $x$ and $u$.
Therefore, formula (\ref{n.13}) implies
\begin{equation}\label{n.14}
\lim_{\de \to 0}\|\de^{-1}\Big(T_{r+\de}(x,
r+\de)-T_{r+\de}(x,r)\Big)+S(r)^{-1}\frac{\partial}{\partial
u}s(x,u)\big|_{u=r}B\wt D^{-1}\|_{L^2}=0.
\end{equation}
Here $\|X(x)\|_{L^2}$ denotes the maximum of the norms in $L^2_m(0,r)$
of the
columns of the matrix function $X$.
Take into account that $T=T^*$, and thus $T_r(x,u)=T_r(u,x)^*$. So, it follows
from 
(\ref{n.8}) and (\ref{n.14}) that
\begin{equation}\label{n.15}
\lim_{\de \to
0}\de^{-1}\int_0^{r}\Big(T_{r+\de}(r+\de,x)-T_{r+\de}(r,x)\Big)\Phi(x)dx=Y_3.
\end{equation}
In a similar to (\ref{n.13}) way, by (\ref{n.2}) and (\ref{n.12})
we have
\begin{equation}\label{n.16}
T_{r+\de}(x,r)-T_{r}(x,r)=-S(r)^{-1}\int_r^{r+\de}s(x,t)T_{r+\de}(t,r)dt.
\end{equation}
From (\ref{n.16}) it follows that
\begin{equation}\label{n.17}
\lim_{\de \to
0}\de^{-1}\int_0^{r}\Big(T_{r+\de}(r,x)-T_{r}(r,x)\Big)\Phi(x)dx=Y_2.
\end{equation}
According to (\ref{n.10}), (\ref{n.11}), (\ref{n.15}), and
(\ref{n.17}) we get
\begin{equation}\label{n.18}
\lim_{\de \to
0}\de^{-1}\Big(\int_0^{r+\de}T_{r+\de}(r+\de,u)
\Phi(u)du-\int_0^{r}T_{r}(r,u)
\Phi(u)du\Big)=Y_1+Y_2+Y_3.
\end{equation}
Using (\ref{n.3}), in a similar way one obtains
\begin{equation}\label{n.19}
\lim_{\de \to 0}\de^{-1}\Big(\int_0^{r}T_{r}(r,u)
\Phi(u)du-\int_0^{r-\de}T_{r-\de}(r-\de,u)
\Phi(u)du\Big)=Y_1+Y_2+Y_3.
\end{equation}
From (\ref{n.18}) and (\ref{n.19}), the equality (\ref{n.9}) is
immediate.
Finally, by (\ref{n.4}) and (\ref{n.9}) we have (\ref{n.5}).

In view of (\ref{n.5})-(\ref{n.8}), $\b$ is differentiable and the
first relations
in (\ref{1.d12}) are valid. According to Remark \ref{Rk3.5} the second
relations
in  (\ref{1.d12}) are fulfilled too. By (\ref{n.4})  we have
$\b(0)=\wt D^{-\frac{1}{2}}\Phi(0)$ and so, taking into account
(\ref{2.4}),
we get
$\b_{k1}(0)=\big(1+|\Phi_{k2}(0)|^2\big)^{-\frac{1}{2}}\not=0$,
that is,
the inequalities (\ref{1.7}) are true.

{\bf Step 2.} Now, let us show that the matrix function
$\b$, which is constructed using (\ref{2.d6}), 
provides the solution of the inverse problem. In view of Definition
\ref{Wf},
Lemma \ref{Law}, and Remark \ref{Rkw} it will suffice to show that
for $r\leq l$ and  $1\leq k \leq m$ we have
\begin{equation}\label{n.27}
\sup_{\Im\mu<-M/4} \left\| w_A(r,\la) \left[
\begin{array}{c}\vp_k(\mu)\\ 1 \end{array}
\right]\right\|  < \infty, \quad
\la=d_k+\frac{b_k}{2\mu}.
\end{equation}
Taking into account (\ref{n.1})  one can see
that the vector function $\mu^{-1}
\vp\left(\frac{\mu}{2}\right)$  is holomorphic in
the half-plane $\eta= \Im \mu<-\frac{M}{2}$ and belongs to
$L^2_m(-
\infty, \, \infty)$ with respect to $\zeta=\Re \mu$ for
any fixed
$\eta<-\frac{M}{2}$. Therefore, similar to the proof of Theorem
\ref{TmInv1} we apply Theorem V from \cite{PW}
and derive:
\begin{equation}\label{2.29.1}
\mu^{-1}
\vp\left(\frac{\mu}{2}\right)={\mathrm{l.i.m.}}_{a
\to
\infty}\int_0^ae^{-i \mu x}e^{\frac{1}{2}x M}f(x)dx, \quad
f \in L^2_m(0,
\infty).
\end{equation}
Using Plansherel's theorem and formula (\ref{2.17}),
 we express $f(x)$  in (\ref{2.29.1})
 via $\Phi_{2}(x)$  and obtain:
\begin{equation}\label{2.31.1}
\frac{i}{\mu}
\vp\left(\frac{\mu}{2}\right)=\int_0^{\infty}e^{-i
\mu x}\Phi_{2}(x)dx, \quad e^{-\frac{1}{2}x M}\Phi_{2}(x)
\in L^2_m(0,
\infty).
\end{equation}
It follows also that the right-hand side of (\ref{2.17'}) equals $0$
for $x<0$,
as stated in the theorem. 
The equality in (\ref{2.31.1}) is true pointwise. 
Recall that according to  (\ref{n.1}) and (\ref{2.17}), $\Phi_2$ 
is differentiable and $e^{\eta u}\Phi_2^{\prime}(u) \in
L^2_m(0, \, \infty)$
for $\eta<-\frac{M}{2}$.
Hence, we rewrite (\ref{2.31.1}) as:
\begin{equation}\label{n.23}
\vp\left({\mu}\right)=-\Phi_2(0)-\int_0^{\infty}e^{-2i
\mu u}\Phi^{\prime}_{2}(u)du. 
\end{equation}
Without loss of generality we shall choose $M$ for the
half-plane $\eta<-M/4$, where $\vp(\mu)$ is treated, so that
\begin{equation}\label{n.25}
\Big(\exp\big(\ve-\frac{M}{2}\big)
u\Big)\Phi_2^{\prime}(u) \in L^2_m(0, \, \infty),
\quad M>4\max_{j\not=p}|d_p-d_j|^{-1}, 
\end{equation}
\begin{equation}\label{n.25'}
 \sup_{\eta\leq-M/4}\|\vp(\mu)\|<\infty.
\end{equation}
In view of (\ref{2.5'}), (\ref{n.21}), and (\ref{n.23}) we have
\begin{equation}\label{n.24}
\big(A_k(r)-\la I\big)^{-1}\Pi_k\left[
\begin{array}{c}\vp_k(\mu)\\ 1 \end{array}
\right]=2b_k\mu e^{2i\mu x}\int_x^\infty e^{-2i \mu
u}\Phi_{k2}^{\prime}(u)du.
\end{equation}
By (\ref{n.24}) and the first relation in (\ref{n.25}) we get
\begin{equation}\label{n.26}
\left\|\big(A_k(r)-\la I\big)^{-1}\Pi_k\left[
\begin{array}{c}\vp_k(\mu)\\ 1 \end{array}
\right]\right\|\leq C(r)|\mu|/\sqrt{|\eta|}.
\end{equation}
By the second relation in (\ref{n.25}) the inequality
$|\la -d_p|>|d_k-d_p|/2$
is true for $p\not=k$, $\la =d_k+\frac{b_k}{2\mu}$, and
$\eta<-M/4$, that is, the resolvents
$\big(A_p(r)-\la I\big)^{-1}$ ($p\not=k$) are bounded.
Therefore, the inequalities (\ref{n.25'})
and  (\ref{n.26})
imply
\begin{equation}\label{n.28}
\left\|\big(A(r)-\la I\big)^{-1}\Pi\left[
\begin{array}{c}\vp_k(\mu)\\ 1 \end{array}
\right]\right\|\leq C_1(r)|\mu|/\sqrt{|\eta|}.
\end{equation}
As $|\la-\ov \la|=-\eta|\mu|^{-2}$,
from (\ref{n.20}), (\ref{n.25'}),  and (\ref{n.28}) the
inequality 
\[
\sup [\ov{ \vp_k(\mu)} \quad 1]w_A(r,\la)^*w_A(r,\la)
\left[
\begin{array}{c}\vp_k(\mu)\\ 1 \end{array}
\right]<\infty
\]
is immediate. Thus, (\ref{n.27}) is valid, and   the solution of the
inverse problem can be obtained via (\ref{2.d6}).

The uniqueness of the solution of the inverse
problem is stated in Theorem \ref{Uniq2}.
\end{proof}
In a way similar to \cite{SaA21', SaA22}, the procedure to solve
the inverse problem grants also a Borg-Marchenko-type
theorem.
\begin{Tm} \label{BoMa}
         Let the  $\BC^m$-valued vector
functions $\varphi(\mu,1)$ and  $\varphi(\mu,2)$ be holomorphic 
in the half-plane
$\Im \mu<-M/4$ $(M>0)$ and
satisfy
(\ref{n.1}). Suppose that on some ray $c \Im \mu =  \Re \mu
$ $\,( c =\ov c$, $ \Im \mu <-M/4)$ we have
\begin{equation} \label{n.30}
\|\varphi(\mu,1)-\varphi(\mu,2)\|=e^{-2i\mu l}O(1)
\quad
{\mathrm{for}} \, |\mu| \to \infty .
\end{equation}

Then $\varphi(\mu,1)$ and  $\varphi(\mu,2)$ are Weyl functions
of systems
(\ref{0.2}) with potentials $\b(x,1)$ and $\b(x,2)$, respectively,
which satisfy  (\ref{1.d12}), (\ref{1.7}) and the additional equality
\begin{equation} \label{n.31}
\b(x,1) \equiv \b(x,2) \quad (0<x<l).
\end{equation}
\end{Tm}
\begin{proof} 
According to Theorem \ref{TmInv2} the functions
$\varphi(\mu,1)$ and  $\varphi(\mu,2)$ are Weyl functions of
systems
(\ref{0.2}) with potentials $\b(x,1)$ and $\b(x,2)$, 
which satisfy  (\ref{1.d12}) and (\ref{1.7}).
Denote by $\Phi_2(x,j)$ the matrix function generated via (\ref{2.17})
by $\varphi(\mu,j)$ $(j=1,2)$, and introduce the matrix function
\begin{equation} \label{n.29}
\nu(\mu):=\int_0^{l}\big(\exp2i\mu(l-u))\big)\Big(
\Phi_2(u,1)-\Phi_2(u,2)\Big) du.
\end{equation}
Next, we shall show that $\nu(\mu)=0$. Indeed,
in view of  (\ref{n.23}) and (\ref{n.30}) it is immediate that
$\Phi_2(0,1)=\Phi_2(0,2)$. Hence, from (\ref{n.23}) and
(\ref{n.29})
we derive
\begin{equation} \label{n.32}
\nu(\mu)=e^{2i\mu l}\big(
\varphi(\mu,2)-\varphi(\mu,1) \big)
+\int_l^{\infty}\big(\exp2i\mu(l-u)\big)\Big(
\Phi_2(u,1)-\Phi_2(u,2)\Big) du.
\end{equation}
It is  immediate from
(\ref{n.29})
that $\|\nu(\mu)\|$ is bounded on the line $\Im \mu
=-M/4$.
Using (\ref{n.30}) and (\ref{n.32}) we see that 
$\|\nu(\mu)\|$
is bounded on the  ray $c \Im \mu =  \Re \mu$ $ (\Im
\mu <-M/4)$.   Thus, by the
Phragmen-Lindel\"of theorem $\|\nu(\mu)\|$ is bounded
in the half-plane
$ \Im \mu \leq -M/4$. Moreover, by (\ref{n.29})
$\|\nu(\mu)\|$ is bounded
for $ \Im \mu > -M/4$, that is, $\nu(\mu)$ is an entire
function bounded
on $\BC$. Therefore,  $\nu(\mu)$ is a constant. As
$\lim_{\eta \to \infty}
\nu(\mu)=0$, so we obtain $\nu(\mu) \equiv 0$, that is
\begin{equation} \label{n.33}
\Phi_2(x,1) \equiv \Phi_2(x,2)\quad (0<x<l).
\end{equation}
By the procedure to solve the inverse problem (see Theorem \ref{TmInv2})
formula (\ref{n.33}) implies (\ref{n.31}).
\end{proof}
\section{Weyl set} \label{Set}
\setcounter{equation}{0}
Condition (\ref{1.7})  is not necessary in the considerations
of Sections \ref{prel},   \ref{S} and  \ref{Inv}.
Taking into account the identity in  (\ref{0.3}), one can choose
two sets  of natural numbers $N_1$ and $N_2$, so that
\begin{equation}\label{y25}
\beta_{k1}(0)\not=0 \quad {\mathrm{for}}  \quad k \in
N_1, \quad 
\beta_{k2}(0)\not=0 \quad {\mathrm{for}}  \quad k \in
N_2,
\end{equation}
\begin{equation}\label{y26}
N_1\cap N_2=\emptyset, \quad N_1\cup N_2=\{1, \, 2,
\, \ldots, \, m  \}.
\end{equation}
The procedure to solve the inverse problem is easily modified for that
case. Such a modification is discussed below. We  also introduce
a notion of the Weyl set, and the system (\ref{0.2}), which satisfies
(\ref{0.3}),
is recovered from its  Weyl set.
\begin{Dn}\label{WS} The set $\vk$ of the pairs
\begin{equation}\label{y27}
\vk=\{\b_k(0), \, \wt \psi_k(\mu)\, | \,
1\leq k \leq m \},
\end{equation}
where $\wt \psi_k$ is given by   (\ref{1.d6}), is called a Weyl set
of the system (\ref{0.2}), which satisfies (\ref{0.3}).
\end{Dn}
If the functions $\wt \psi_k$ in the Weyl sets $\vk$ and
$\breve \vk$
coincide and the vectors $\b_k(0)$ and $\breve\b_k(0)$ in 
$\vk$ and $\breve \vk$, respectively, satisfy the equalities
$\breve\b_k(0)=c_k\b_k(0)$, $|c_k|=1$ ($1\leq k \leq m$),
we say that $\vk=\breve \vk$, as $\vk$ and $\breve \vk$
correspond to the same system.
After we exclude this arbitrariness, the Weyl set is uniquely defined
by system  (\ref{0.2}),  (\ref{0.3}). (See the construction of $\wt
\psi_k$
in Section \ref{prel} and the proof of Theorem  2.2 in \cite{MST}.)

The next lemma is proved in a quite similar way to Lemma \ref{Law}.
\begin{La}\label{Laws} Let the $S$-node be given by the formulas
(\ref{2.1})-(\ref{2.3}) and (\ref{2.5}), where
\begin{equation}\label{y28}
\Phi_{k1}(x)\equiv 1, \quad  \Phi_{k2}^{\prime}(x)\in
L^2(0,l) \quad {\mathrm{for}}  \quad k \in N_1,
\end{equation}
\begin{equation}\label{y29}
\Phi_{k2}(x)\equiv 1, \quad  \Phi_{k1}^{\prime}(x)\in
L^2(0,l) \quad {\mathrm{for}}  \quad k \in N_2.
\end{equation}
Define the reductions $A(r)$, $S(r)$ and $\Pi(r)$ $(r<l)$ of the
operators from this $S$-node
in the same way
as it was done after Remark \ref{res} in Section \ref{S}.

Then $w_A$ of the form (\ref{2.d4}) satisfies (\ref{2.d5}), where
$\b$ is given
by  (\ref{2.d6}) and the operator $V$ in (\ref{2.d6}) is obtained via 
(\ref{2.13})-(\ref{2.15}).
\end{La}
\begin{Rk}\label{RkOI}
The operator $S$ is uniquely recovered from the operator identity
(\ref{2.5}),
the corresponding formulas from Proposition \ref{Pn2.1} are easily modified
for the more general case  (\ref{y28}), (\ref{y29}).
\end{Rk}
In view of  Lemma \ref{Laws}, the proof of the next theorem
is similar to the proof of Theorem \ref{TmInv1}.
\begin{Tm}\label{TmInv3}
Let system (\ref{0.2}) satisfy conditions (\ref{0.3}) and
let $N_1$ and $N_2$ be chosen so that  (\ref{y25}) and (\ref{y26})
hold.

Then, system (\ref{0.2}) is uniquely recovered from its Weyl set.
To
recover (\ref{0.2}) on an arbitrary fixed interval $[0, \, l]$ we
use  the following procedure.

First, the entries of the column vector functions $\Phi_1(x)$
and  $\Phi_2(x)$
$(0<x<\infty)$  are introduced by the equalities:
\begin{equation}\label{y30}
\Phi_{k1}(x)\equiv 1, \quad
\Phi_{k2}(x)=\frac{i}{2
\pi}\int_{-\infty}^{\infty}\mu^{-1}e^{i \mu
x}\vp_k\left(\frac{\mu}{2}\right)d\zeta  \quad
{\mathrm{for}}  \, k \in N_1;
\end{equation}
\begin{equation}\label{y31}
\Phi_{k1}(x)=\frac{i}{2
\pi}\int_{-\infty}^{\infty}\mu^{-1}e^{i \mu
x}\phi_k\left(\frac{\mu}{2}\right)d\zeta,   \quad
\Phi_{k2}(x)\equiv 1 \quad {\mathrm{for}}  \, k \in N_2;
\end{equation}
\begin{equation}\label{y32}
\vp_k(\mu)=\frac{\ov{\b_{k1}(0)}\wt
\psi_k(\mu)-\b_{k2}(0)}{\ov{\b_{k2}(0)}\wt
\psi_k(\mu)+\b_{k1}(0)},
 \quad  \phi_k(\mu)=\frac{\ov{\b_{k2}(0)}\wt
\psi_k(\mu)+\b_{k1}(0)}{\ov{\b_{k1}(0)}\wt
\psi_k(\mu)-\b_{k2}(0)},
\end{equation}
where $\mu = \zeta +i \eta$ and $-\eta>0$ is sufficiently
large.

Next, we introduce the $S$-node by formulas (\ref{2.1})-(\ref{2.3}) and
(\ref{2.5}).
The operator $S$ has the form (\ref{2.6}), where the definition
(\ref{2.6'}) of $\wt D$
is modified:
\begin{equation}\label{y33}
\wt D=\wt D_1+\wt D_2, \quad \wt D_p
={\rm{diag}}\{|\Phi_{1p}(0)|^2,|\Phi_{2p}(0)|^2,\ldots,|\Phi_{mp}(0)|^2\}
\quad (p=1,2).
\end{equation}

Finally, we put
\begin{equation}\label{y34}
\breve \b(x)=\left[
\begin{array}{c}
\breve \b_1(x) \\ \cdots \\ \breve \b_m(x)
\end{array}
\right]=(V\Phi)(x)  \quad (0 \leq x \leq l),
\end{equation}
where the operator $V$  is obtained via  (\ref{2.13})-(\ref{2.15}).

The equalities
\begin{equation}\label{y35}
\breve \b_k(x)^* \breve \b_k(x) \equiv 
 \b_k(x)^*  \b_k(x) \quad (1\leq k \leq m).
\end{equation}
are valid, that is, system   (\ref{0.2}) is recovered by the procedure, which
is
described above.
\end{Tm}

\section{Sine-Gordon equation} \label{Sine}
\setcounter{equation}{0}
The initial value problem for the sine-Gordon equation in the light cone coordinates $\om_{xt}=\sin
\, \om$ ($\om_x:=\frac{\p }{\p x}\om$)
was treated in \cite{AKNS}. The initial value problem (with initial
conditions tending to zero)
for the sine-Gordon equation in laboratory coordinates
\begin{equation}\label{y36}
\omega_{xx} - \omega_{tt} = \sin \omega
\end{equation}
was investigated by  Faddeev, Takhtajan and Zakharov  (see \cite{ZTF} and further
references in \cite{TF}).  Notice also that the Goursat problem for the equation $\om_{xt}=\sin
\, \om$, which is treated on the characteristics $t=0$ and $\,x=-\infty$ in \cite{KN1}, is equivalent
to the Cauchy problem for equation \eqref{y36}. In this section we consider  \eqref{y36}
under boundary conditions $\om(0,t)=\om_0(t)$ and  $ \om_x(0,t)=\om_1(t)$. (The boundary value
problem is clearly equivalent to the initial value problem after the change of variables and the change
$\om \to \om+\pi$.) We do not require that $\om $ tends to zero and only the boundedness
of $\om_x$ and $\om_t$ is needed.

Equation (\ref{y36}) admits  a  zero curvature representation
\begin{equation} \label{y37}
G_t(x,t,\la)-F_x(x,t, \la)+G(x,t,\la)
F(x,t,\la)-F(x,t,\la)G(x,t,\la)=0.
\end{equation}
We can modify the auxiliary systems
\begin{equation} \label{y38}
w_x(x,t, \la)=G(x,t, \la)w(x,t,\la), \quad w_t(x,t,
\la)=F(x,t,\la)w(x,t,\la),
\end{equation}
so that they have the form (\ref{0.2}). Namely, put
\begin{eqnarray}\nonumber
&&G(x,t,\la)= i \sum_{k=1}^2 b_k(\la-d_k)^{-1}\Big(
\beta_k(x,t)^* \beta_k(x,t)\Big), \\ && 
d_1=-d_2=1, \quad b_1=b_2=1, \label{y39}
\end{eqnarray}
\begin{eqnarray}\nonumber
&&F(x,t,\la)= i \sum_{k=1}^2 b_k(\la-d_k)^{-1}\Big(
\beta_k(x,t)^* \beta_k(x,t)\Big), \\ && d_1=-d_2=1,
\quad b_1=-b_2=1, \label{y40}
\end{eqnarray}
where
\begin{equation}\label{y41}
\beta_1(x,t) = \frac{ 1 }{\sqrt 2}[1 \quad i \, e^{i
\omega(x,t)/2}]
q(x,t), \qquad
\beta_2(x,t) = \frac{ 1} {\sqrt 2} [1 \quad i e^{-i \,
\omega(x,t)/2}]
q(x,t),
\end{equation}
the $2 \times 2$ matrix function $q$ satisfies the equations
\begin{equation}\label{y42} 
q_x(x,t)=\breve G(x,t) q(x,t), \quad  q_t(x,t)=\breve F(x,t)
q(x,t), \quad q(0,0)=I_2,
\end{equation}
\begin{equation}\label{y43}
\breve G:= - i \left( \displaystyle \frac{\omega_t
\,}{4} j + \frac{ 1}{2} \sin
\left( \frac {\omega}
{2} \right) J \right), \quad
\breve F: =  - i \displaystyle \frac{\omega_x}{4} j +
\frac{ 1}{2 }\cos
\left( \frac{ \omega}
{2} \right) J j,
\end{equation}
and $J$, $j$ are given in (\ref{y16'}). It is easily checked that the
sine-Gordon
equation  (\ref{y36}) is equivalent to the compatibility condition
$\breve G_t-\breve F_x+\breve G\breve F-\breve F\breve
G=0$
of  the equations  (\ref{y42}). Moreover, direct calculation shows that  
relations (\ref{y41}) and  (\ref{y42}) imply  (\ref{y37}),
which is the compatibility condition for  (\ref{y38}). Thus, if 
(\ref{y36}) holds, equations  (\ref{y38}) are compatible.

Introduce the $2 \times 2$ matrix functions $Z(x,t,\la)$,
\[ Z(t,\la):=Z(0,t,\la)=\{Z_{ij}(t,\la)\}_{i,j=1}^2, \quad
{\mathrm{and}} \quad
Y(x,t,\la)=\{Y_{ij}(x,t,\la)\}_{i,j=1}^2\]
 by the
equations
\begin{eqnarray}\nonumber &&
Y_x(x,t,\la)=G(x,t,\la)Y(x,t,\la), \quad Y(0,t,\la)\equiv
I_2; \\ &&Z_t(x,t,\la)=F(x,t,\la)Z(x,t,\la), \quad
Z(x,0, \la)\equiv I_2.
\label{y44}
\end{eqnarray}
The matrix functions $Q_k(x,t)$  ($k=1,2$) are  connected with $\b_k(x,t)$
by the equalities (\ref{1.2}).
According to  (\ref{1.2}) and  (\ref{y41})-(\ref{y43}) the boundary
conditions
\begin{equation}\label{y45}
\om(0,t)=\om_0(t), \quad  \om_x(0,t)=\om_1(t) \quad
(-\infty<t<\infty)
\end{equation}
uniquely define $\breve F(0,t)$, $q(0,t)$, $\b_k(0,t)$  and $Q_k(0,t)$.
If we recover also $\wt \psi_k(t, \mu)$ for each
$-\infty<t<\infty$, we  have
a Weyl set for each $t$. 
First, we recover $F(0,t, \la)$ and $Z(t,\la)$, using formulas 
(\ref{y40}) and (\ref{y44}), and
put
\begin{equation}\label{y46}
U_k(x,t,\mu):=\exp\{(-1)^kit
\mu\}Q_k(x,t)Z(x,t,\la)Q(x,0)^*, \quad
\mu=\big(2(\la -d_k)\big)^{-1},
\end{equation}
\begin{equation}\label{y47}
U_k(t,\mu)=\{u_{jp}(t,\mu,k)\}_{j,p=1}^2:=U_k(0,t,\mu).
\end{equation}
The matrix functions $U_k(t,\mu)$ are uniquely  recovered from
(\ref{y45}) too.
\begin{Tm}\label{TmSG} Let the function $\om(x,t)$ have continuous
second derivatives
in the semi-plane $x\geq 0$ and satisfy  the sine-Gordon equation
(\ref{y36}) and boundary conditions (\ref{y45}). 
Assume  also that
\begin{equation}\label{y48}
\sup_{x\geq 0}(|\om_x(x,t)|+|\om_t(x,t))|)<\infty .
\end{equation}

Then,  $\cos \om(x,t)$ $(x\geq 0)$   is uniquely  recovered
from  (\ref{y45}).
For this purpose
construct $U_k(t,\mu)$ $(k=1,2)$ using (\ref{y40})-(\ref{y47}). 
There is $M_1>0$,
such that for
$- \Im \mu >M_1$ we have
\begin{equation}\label{y49}
\wt \psi_1(0, \mu)=-\lim_{t\to \infty}
\frac{u_{12}(t,\mu, 1)}{u_{11}(t,\mu,1)}, \quad 
\wt \psi_2(0, \mu)=-\lim_{t\to -\infty}
\frac{u_{12}(t,\mu, 2)}{u_{11}(t,\mu,2)}.
\end{equation}
The functions  $\wt \psi_k(t, \mu)$ are given by the formulas
\begin{equation}\label{y50}
\wt \psi_k(t, \mu)= \frac{u_{11}(t,\mu, k)\wt
\psi_k(0, \mu)+u_{12}(t,\mu, k) }{u_{21}(t,\mu, k)\wt
\psi_k(0, \mu)+u_{22}(t,\mu,k)}.
\end{equation}
By formulas  (\ref{1.2}),  (\ref{y41})-(\ref{y43}), 
(\ref{y49}) and  (\ref{y50}) we recover
the Weyl set for each $t$. 
Finally, we recover the functions $\b_k(x,t)$ $($up to  factors $c_k(x,t)$
such that $|c_k|=1$, $c_k(0,t)=1$ $)$ using Theorem \ref{TmInv3}.
It follows that
\begin{equation}\label{y51}
\cos
\om(x,t)=2\b_1(x,t)\b_2(x,t)^*\b_2(x,t)\b_1(x,t)^*-1.
\end{equation}
\end{Tm}
\begin{proof}  
{\bf Step 1.}  In this step we shall prove  (\ref{y50}).
Note that as $\om$ has continuous second derivatives, so according to
(\ref{y39})-(\ref{y43})
the matrix functions $G$ and $F$ are continuously differentiable. Therefore,
the
formula (1.6) on p. 168 in  \cite{SaL3} implies:
\begin{equation}\label{y52}
Y(x,t,\la)=Z(x,t,\la)Y(x,0,\la)Z(t,\la)^{-1}.
\end{equation}
By (\ref{d1}),  (\ref{1.3}), and (\ref{y44}) we have $w(x,t,
\la)=Y(x,t,\la)$,
\begin{equation}\label{y53}
 Q(x,t)^*=Q(x,t)^{-1},   \quad Y(x,t,\ov \la)^*=Y(x,t, \la)^{-1},
\quad Z(x,t,\ov \la)^*=Z(x,t, \la)^{-1}.
\end{equation}
Hence, taking into account (\ref{1.0}) and  (\ref{1.d1}) we have
\begin{equation}\label{y54}
{\mathfrak A}(r,t, \mu)=e^{ir\mu}Q(0,t) Y(r,t,\ov
\la)^*Q(r,t)^*.
\end{equation}
From (\ref{y46}), (\ref{y47}), (\ref{y52}), and (\ref{y54}) it
follows that
\begin{equation}\label{y55}
{\mathfrak A}_k(r,t, \mu)=U_k(t,\mu){\mathfrak A}_k(r,0,
\mu)U_k(r,t,\mu)^{-1}.
\end{equation}
In view of  (\ref{y44}) and (\ref{y46}) it is easy to see that
\begin{eqnarray}\label{y56}
&&\frac{\p}{\p
t}U_k(r,t,\mu)=\Big((-1)^{k+1}i\mu
j+\Big(\frac{\p}{\p t}Q_k(r,t)\Big)Q_k(r,t)^*
\\ \nonumber &&
+
(-1)^{k}iQ_k(r,t)
\frac{\b_p(r,t)^*\b_p(r,t)}{\la -d_p}Q_k(r,t)^*\Big)U_k(r,t,\mu),
\end{eqnarray}
where  $k$ and $p$ take values $1$ and $2$, $p\not=k$.  
According to  (\ref{y42}) and  (\ref{y43}) the equalities
$ \frac{\p}{\p t}(q^*q)=0$ and $ \frac{\p}{\p
x}(q^*q)=0$ are true
and $q(0,0)=I_2$. Therefore, it is immediate that $q$ is unitary:
\begin{equation}\label{y56'}
q(x,t)^*q(x,t) \equiv I_2.
\end{equation}
By  
 (\ref{y41})-(\ref{y43}),   (\ref{y48}) and  (\ref{y56'}) 
we have
\begin{equation}\label{y57}
\sup \|\frac{\p}{\p t}\b_k(x,t)\|<\infty
\quad (x\geq 0, \, -\infty<t<\infty),
\quad  k=1,\, 2.
\end{equation}
Taking into account
(\ref{y56}) and (\ref{y57}), in a way  similar to the proof of  (\ref{1.10}) we derive
\begin{equation}\label{y58}
(-1)^{k+1}\frac{\p}{\p
t}\Big(U_k(r,t,\mu)^*jU_k(r,t,\mu)\Big)\geq \big(i(\mu -
\ov \mu) -\wh M\big)U_k(r,t,\mu)^*U_k(r,t,\mu)>0
\end{equation}
for some $\wh M>0$ and $\Im \mu<-\wh M/2$. From  
(\ref{y58}) it follows that 
\begin{equation}\label{y58'}
U_1(r,t,\mu)^*jU_1(r,t,\mu) > j \quad {\mathrm{for}} \,
t> 0; \quad U_2(r,t,\mu)^*jU_2(r,t,\mu) > j \quad
{\mathrm{for}} \, t< 0.
\end{equation}
From the first inequality in \eqref{y58'} we have
$\big(U_1(r,t,\mu)^*\big)^{-1}jU_1(r,t,\mu)^{-1} < j $,
\begin{equation}\label{y59}
 [\t(\mu)^* \quad 1]
 \big(U_1(r,t,\mu)^*\big)^{-1}jU_1(r,t,\mu)^{-1}\left[
\begin{array}{c}
\t(\mu) \\ 1
\end{array}
\right]<0 \quad
{\mathrm{for}} \, \, t> 0, \, |\t|\leq 1. 
\end{equation}
In a similar way, from the second inequality in  (\ref{y58'}) we derive
\begin{equation}\label{y60}
[\t(\mu)^* \quad 1]
 \big(U_2(r,t,\mu)^*\big)^{-1}jU_2(r,t,\mu)^{-1}\left[
\begin{array}{c}
\t(\mu) \\ 1
\end{array}
\right]<0 \quad
{\mathrm{for}} \, \, t< 0, \, |\t|\leq 1. 
\end{equation}
By  (\ref{y59})  and   (\ref{y60}) we obtain 
the inequalities 
\begin{equation}\label{y61}
|\chi_1(r,t,\mu)|<1 \quad 
{\mathrm{for}} \quad  t> 0; \quad |\chi_2(r,t,\mu)|<1
\quad 
{\mathrm{for}} \quad  t< 0 \quad ( \Im \mu
<-\frac{\wh M}{2})
\end{equation}
for  the functions
\[
\chi_k(r,t,\mu):=\frac{(U_k^{-1})_{11}(r,t,\mu)\t(\mu)+
(U_k^{-1})_{12}(r,t,\mu)}{(U_k^{-1})_{21}(r,t,\mu)
\t(\mu)+(U_k^{-1})_{22}(r,t,\mu)},  \quad k=1,2; \,\,
|\t|\leq 1.
\]
In view of  (\ref{1.d2})  and   (\ref{1.d6})  we have
\begin{equation}\label{y62}
\wt \psi_k(t,\mu)=\lim_{r\to
\infty}\frac{{\mathfrak A}_{11}(r,t,\mu)\t(\mu)+
{\mathfrak
A}_{12}(r,t,\mu)}{{\mathfrak A}_{21}(r,t,\mu)
\t(\mu)+{\mathfrak
A}_{22}(r,t,\mu)}, \quad |\t(\mu)|\leq 1,  \quad
\Im \mu
<-\frac{M}{4}.
\end{equation}
From (\ref{y55}) it follows that the linear fractional transformation on
the right-hand side
of  (\ref{y62}) can be written down as the superposition of the three
linear fractional transformations, the first of which transforms  $\t$ into
$\chi_k$.
Using (\ref{y61}) we see that the second transformation transforms
$\chi_k$
into $\psi_k(r,0, \mu)$ for $k=1$ and $t>0$ as well as for  $k=2$
and $t<0$.
In the limit, it follows from (\ref{y55}), (\ref{y61}), and
(\ref{y62}) that
(\ref{y50}) is true  for $k=1$, $t>0$ and for $k=2$, $t<0$
$\big(-\eta>M_0=\max(\frac{\wh M}{2}, \, \frac{
M}{4})\big)$.

To prove (\ref{y50})   for $k=1$, $t<0$ and for $k=2$, $t>0$ rewrite
(\ref{y55}) in the form
\begin{equation}\label{y63}
{\mathfrak A}_k(r,t, \mu)U_k(r,t,\mu)=U_k(t,\mu){\mathfrak
A}_k(r,0, \mu),
\end{equation}
and use the inequalities
\begin{equation}\label{y64}
U_1(r,t,\mu)^*jU_1(r,t,\mu) < j \, \, {\mathrm{for}}
\,\,  t< 0; \quad U_2(r,t,\mu)^*jU_2(r,t,\mu) < j 
\, \, {\mathrm{for}} \, \, t> 0,
\end{equation}
which are immediate from  (\ref{y58}). By  (\ref{y64}) one can see that
\begin{equation}\label{y65}
|\breve \chi_1(r,t,\mu)|<1 \quad 
{\mathrm{for}} \quad  t< 0; \quad |\breve
\chi_2(r,t,\mu)|<1 \quad 
{\mathrm{for}} \quad  t> 0 \quad (  \Im \mu
<-\frac{\wh M}{2}),
\end{equation}
where
\[
\breve
\chi_k(r,t,\mu):=\frac{(U_k)_{11}(r,t,\mu)\t(\mu)+
(U_k)_{12}(r,t,\mu)}{(U_k)_{21}(r,t,\mu)
\t(\mu)+(U_k)_{22}(r,t,\mu)},  \quad k=1,2; \,\,
|\t|\leq 1.
\]
Now, consider the linear fractional transformations of $\t$, where the
coefficients are the entries of the
left-hand side and right-hand side of  (\ref{y63}), respectively. These
linear fractional transformations coinside,
and in the limit (as $r$ tends to infinity) we obtain (\ref{y50})   for
$k=1$, $t<0$ and for $k=2$, $t>0$.
Thus, (\ref{y50})  is proved.

{\bf Step 2.}  By (\ref{1.d3}) we have $|\wt
\psi_k(t,\mu)|<1$  ($\Im \mu <-\frac{ M}{4}$).  Hence,
in view of  (\ref{y50})  we obtain
\begin{equation}\label{y66}
[\wt \psi_k(0,\mu)^* \quad 1]
 U_k(t,\mu)^*jU_k(t,\mu)\left[
\begin{array}{c}
\wt \psi_k(0,\mu)\\ 1
\end{array}
\right] \leq 0. 
\end{equation}
Recall that (\ref{y50})  holds for $-\Im \mu>M_0 \geq
\frac{\wh M}{2}$.
From  (\ref{y58}) it follows that for $-\Im
\mu>M_1=\frac{\ve}{2}+M_0$ ($\ve>0$) the
inequalities
\begin{eqnarray} \label{y67!}
 &&U_1(t,\mu)^*jU_1(t,\mu)-j \geq \ve \int_0^t
U_1(s,\mu)^*U_1(s,\mu)ds \quad (t>0),
 \\
 &&\label{y67}
 U_2(t,\mu)^*jU_2(t,\mu)-j \geq \ve \int_t^0
U_2(s,\mu)^*U_2(s,\mu)ds \quad (t<0)
 \end{eqnarray}
 are true.
According to  (\ref{y66})--(\ref{y67}),  for any $\mu$ such that
$-\Im \mu>M_1$ we have
\begin{eqnarray} \label{y68!}
&&\int_0^\infty[\wt \psi_1(0,\mu)^* \quad 1]
U_1(s,\mu)^*U_1(s,\mu)\left[
\begin{array}{c}
\wt \psi_1(0,\mu)\\ 1
\end{array}
\right] ds<\infty, \\ &&
\int_{-\infty}^0 [\wt \psi_2(0,\mu)^* \quad 1]
U_2(s,\mu)^*U_2(s,\mu)\left[
\begin{array}{c}
\wt \psi_2(0,\mu)\\ 1
\end{array}
\right] ds<\infty.  \label{y68}
\end{eqnarray}
By  (\ref{y48}), (\ref{y56}), (\ref{y68!}),  and (\ref{y68}) the inequalities
\begin{equation}\label{y69}
\sup_{t>0}\left\|U_1(t,\mu)\left[
\begin{array}{c}
\wt \psi_1(0,\mu)\\ 1
\end{array}
\right]   \right\|<\infty, \quad
\sup_{t<0}\left\|U_2(t,\mu)\left[
\begin{array}{c}
\wt \psi_2(0,\mu)\\ 1
\end{array}
\right]   \right\|<\infty
\end{equation}
are valid. Inequalities  (\ref{y67!}) and (\ref{y67}) imply that
\begin{equation}\label{y70}
|u_{11}(t,\mu,1)|^2>1+\ve t \quad (t>0), \quad
|u_{11}(t,\mu,2)|^2>1-\ve t \quad (t<0).
\end{equation}
From  (\ref{y69})  and (\ref{y70})  follows (\ref{y49}).

In view of  (\ref{y41})  and (\ref{y56'})  we obtain 
\[
2\b_1 \b_2^*=1+\cos \om + i \sin \om, \quad
{\mathrm{i.e.,}} \quad
2|\b_1 \b_2^*|^2=1+\cos \om.
\]
Hence, the equality  (\ref{y51}) is immediate.
\end{proof}

Fakult\"at f\"ur Mathematik, Universit\"at Wien,
\\
Nordbergstrasse 15, A-1090 Wien,  Austria. \\
al$_-$sakhnov@yahoo.com

\end{document}